\documentclass[a4paper,12pt,twoside]{article}
\usepackage{graphicx}\usepackage{a4wide}
\usepackage{amssymb,amsmath,amsthm}
\usepackage{subfigure,color,colordvi,graphicx,xcolor}
\usepackage[only,llbracket,rrbracket,interleave]{stmaryrd}
\usepackage[sort,compress]{cite}
\usepackage[pdftex,colorlinks,bookmarksnumbered,bookmarksopen,citecolor=red,urlcolor=red]{hyperref}
\usepackage{multirow}

\newcommand{\bm}[1]{\text{\boldmath $#1$\unboldmath}}

\newcommand{\curl}{\operatorname{curl}}

\newcommand{\mat}[1]{\mathbf{#1}}

\newcommand{\grad}{\bm{\nabla}}
\newcommand{\RR}{\mathbb{R}}

\newcommand{\VhHat}{\ensuremath{\mathcal{\hat{V}}^h}}
\newcommand{\Vh}{\ensuremath{\mathcal{V}^h}}

\newcommand{\eltwo}{\ensuremath{\mathcal{L}_2}}

\newcommand{\hone}{\ensuremath{\mathcal{H}^1}}

\newcommand{\ndof}  {\ensuremath{\texttt{n}_{\texttt{dof}}}}

\newcommand{\nsd}  {\ensuremath{\texttt{n}_{\texttt{sd}}}}
\newcommand{\msd}  {\ensuremath{\texttt{m}_{\texttt{sd}}}}
\newcommand{\numel}{\ensuremath{\texttt{n}_{\texttt{el}}}}

\DeclareMathOperator{\tr}{tr}

\newcommand{\hv}{\hat{v}}

\newcommand{\bu}{\bm{u}}

\newcommand{\bt}{\bm{t}}

\newcommand{\bhw}{\widehat{\bw}}
\newcommand{\bhu}{\widehat{\bu}}

\newcommand{\bv}{\bm{v}}

\newcommand{\bn}{\bm{n}}
\newcommand{\bx}{\bm{x}}
\newcommand{\bL}{\bm{L}}

\newcommand{\btau}{\bm{\tau}}
\newcommand{\btauE}{\bm{\tau}_{\!e}}

\newcommand{\nDeg}{\ensuremath{k}}
\newcommand{\Pk}{\ensuremath{\mathcal{P}^{\nDeg}}}
\newcommand{\Hdiv}{H(\operatorname{div})}

\newcommand{\Insd}{\mat{I}_{\nsd}}

\newcommand{\stress}{\bm{\sigma}}
\newcommand{\vort}{\bm{\omega}}
\newcommand{\defo}[1]{\bm{\nabla}^{\texttt{S}}#1}
\newcommand{\gradS}{\bm{\nabla}_{\!\texttt{S}}}
\newcommand{\gradW}{\bm{\nabla}_{\!\texttt{W}}}
\newcommand{\stressV}{\bm{\sigma}_{\texttt{V}}}
\newcommand{\strainV}{\bm{e}_{\texttt{V}}}
\newcommand{\bD}{\mat{D}}
\newcommand{\bDHalf}{\bD^{1/2}}
\newcommand{\bE}{\mat{E}}
\newcommand{\bN}{\mat{N}}
\newcommand{\bT}{\mat{T}}

\newcommand{\nrr}  {\ensuremath{\texttt{n}_{\texttt{rr}}}}
\newcommand{\bTens}{\bm{\varsigma}}
\newcommand{\bTensV}{\bm{\varsigma}_{\texttt{V}}}

\newcommand{\jump}[1]{\llbracket #1\rrbracket}

\newcommand{\bw}{\bm{w}}

\newenvironment{keywords}{\begin{quote}\emph{\textbf{Keywords:}}}{\end{quote}}

\newtheorem{lemma}{Lemma}
\theoremstyle{definition}

\newtheorem{remark}{Remark}

\graphicspath{{figs/}}

\begin{document}
\title{A superconvergent HDG method for Stokes flow with strongly enforced symmetry of the stress tensor}

\author{
			  Matteo Giacomini\\[-1ex]
             \small Laboratori de C\`alcul Num\`eric (LaC\`aN), \\[-1ex]
             \small ETS de Ingenieros de Caminos, Canales y Puertos, \\[-1ex]
             \small Universitat Polit\`ecnica de Catalunya, Barcelona, Spain\\[1em]
             Alexandros Karkoulias\\[-1ex]
			  \small Laboratori de C\`alcul Num\`eric (LaC\`aN), \\[-1ex]
             \small ETS de Ingenieros de Caminos, Canales y Puertos, \\[-1ex]
             \small Universitat Polit\`ecnica de Catalunya, Barcelona, Spain\\[-1ex]
             \small Centre Internacional de Metodes Numerics en Enyinyeria (CIMNE), \\[-1ex]
             \small Campus Nord UPC, Barcelona, Spain\\[1em]
			  Ruben Sevilla\\[-1ex]
             \small Zienkiewicz Centre for Computational Engineering, \\[-1ex]
             \small College of Engineering, Swansea University, Wales, UK \\[1em]
             Antonio Huerta\\[-1ex]
             \small Laboratori de C\`alcul Num\`eric (LaC\`aN), \\[-1ex]
             \small ETS de Ingenieros de Caminos, Canales y Puertos, \\[-1ex]
             \small Universitat Polit\`ecnica de Catalunya, Barcelona, Spain\\[1em]
}
\date{\today}
\maketitle

\begin{abstract}
This work proposes a superconvergent hybridizable discontinuous Galerkin (HDG) method for the approximation of the Cauchy formulation of the Stokes equation using same degree of polynomials for the primal and mixed variables.
The novel formulation relies on the well-known Voigt notation to strongly enforce the symmetry of the stress tensor.
The proposed strategy introduces several advantages with respect to the existing HDG formulations.
First, it remedies the suboptimal behavior experienced by the classical HDG method for formulations involving the symmetric part of the gradient of the primal variable. 
The optimal convergence of the mixed variable is retrieved and an element-by-element post-process procedure leads to a superconvergent velocity field, even for low-order approximations.
Second, no additional enrichment of the discrete spaces is required and a gain in computational efficiency follows from reducing the quantity of stored information and the size of the local problems.
Eventually, the novel formulation naturally imposes physical tractions on the Neumann boundary. 
Numerical validation of the optimality of the method and its superconvergent properties is performed in 2D and 3D using meshes of different element types.
\end{abstract}


\section{Introduction and motivations}
\label{sc:Intro}

The interest in discontinuous Galerkin (DG) methods \cite{DGbook2000,MR2372235,Riviere2008,MR2882148,MR3204392,Cangiani2017} has increased in the past years owing to their ability to construct high-order discretizations on unstructured meshes and to their flexibility in performing $p$-adaptivity.
Among the different techniques proposed in the literature to approximate incompressible flow problems, the recent growing interest towards hybridizable discontinuous Galerkin (HDG) methods \cite{Nguyen-NPC:11,MR3626531} is due to multiple advantages these formulations have with respect to classical DG ones.
Concerning Stokes flow, several HDG formulations have been proposed in the literature \cite{Nguyen-NPC:10,Nguyen-CNP:10,MR2904582} and the interested reader is referred to \cite{MR3209969} for an overview on the topic.

The use of hybridization was first introduced with the local discontinuous Galerkin (LDG) method to circumvent the construction of divergence-free approximations of the velocity field \cite{MR2196980}, see also \cite{MR2412047,MR2407142,AdM-MFH:08,MR2683650}.
Moreover, owing to hybridization \cite{Jay-CG:05,Jay-CG:05b}, the globally coupled unknowns are defined on the boundary of the mesh elements and are connected solely to neighboring elements. Thus, the size of the global problem is greatly reduced.
In addition, it is worth noting that HDG allows equal interpolation for velocity, pressure and strain rate tensor, owing to an appropriate definition of the numerical flux and to the introduction of a stabilization parameter.
Thus, the limitations of using equal-order approximations for velocity and pressure in the incompressible limit, through the fulfillment of the Ladyzhenskaya-Babu\v{s}ka-Brezzi (LBB) condition, are circumvented by HDG.
In particular, Cockburn and co-workers~\cite{Cockburn-CDG:08} proved solvability and stability under the aforementioned assumptions, without the need of an enriched space for the mixed variable, or a reduced space for the hybrid one.
In \cite{MR2485446,MR2772094}, optimal convergence rates of order $k+1$ are obtained for all the variables using equal degree of approximation $k$, whereas classical DG display suboptimal convergence of order $k$ for the pressure and the gradient of the velocity.

A key aspect of HDG is the ability to construct a post-processed velocity field superconverging with order $k+2$ \cite{MR3008833}. 
This is crucial when the superconvergent solution is sought to devise automatic procedures to perform $p$-adaptivity (cf. e.g. \cite{giorgiani2014hybridizable,RS-SH:18}).
Nevertheless, the classical HDG equal-order approximation for the Cauchy formulation is known to experience suboptimal convergence of the mixed variable and a loss of superconvergence of the post-processed velocity field using low-order approximations \cite{Nguyen-CNP:10}.

Recently, in a series of publications \cite{MR3626530,MR3601005,MR3601012,MR3649424}, Cockburn and co-workers devoted a great effort to develop a general framework, namely the $\bm{M}$-decomposition, to devise superconvergent HDG discretizations.
This approach relies on enriching the local spaces for the approximation of the mixed variable by adding extra basis functions.
The number of these additional basis functions is not significantly big and in most cases it depends on the type of element under analysis and not on the degree of approximation $k$.
Despite only the size of the local problems increases and the additional computational effort is limited, it induces a more complex implementation compared to standard HDG methods.
Alternative HDG formulations achieve convergence of order $k+2$ for the velocity field when polynomials of degree $k$ are chosen to approximate the hybrid variable \cite{MR3486522,MR3556409,MR3511719}. These methods rely on utilizing smaller spaces for the mixed variable and larger ones for the velocity and exploiting a special stabilization function, the so-called \emph{reduced stabilization}, to handle them.
Closely related approaches, namely the hybrid high-order (HHO) \cite{MR3283758} and the hybridized weak Galerkin (HWG) \cite{Zhai2015} methods can also achieve the same orders of convergence. 

The present work devises a superconvergent hybridizable discontinuous Galerkin method for the Cauchy formulation of the Stokes equation using the same degree of approximation for the primal and mixed variables.
The rest of this paper is organized as follows.
In Section~\ref{sc:StokesSym}, first, the equations governing the Stokes flow are recalled. 
Then, according to the rationale introduced in \cite{preprintVoigtElasticity} for the linear elasticity equation, the symmetry of the stress tensor is strongly enforced by means of a technique well known in the computational mechanics community, namely the Voigt notation for symmetric tensors.
The corresponding Cauchy formulation of the Stokes equation with strongly enforced symmetry of the stress tensor is derived.
In Section~\ref{sc:HDG}, an HDG discretization is introduced.
A local post-process procedure providing a superconvergent velocity field even for low-order approximations is discussed without resorting to the complex framework of the $\bm{M}$-decomposition.
Moreover, contrary to other HDG formulations, the proposed method features a reduced number of degrees of freedom for the mixed variable and is computationally more efficient since the resulting local problems are smaller.
The novel HDG formulation is validated in Section~\ref{sc:simulations}. Extensive analysis of the optimal convergence and superconvergence rates of the primal, mixed and post-processed variables, for two and three dimensional problems is provided by means of numerical simulations. Special emphasis is placed on the influence of the stabilization parameter and on the robustness of the method using meshes of different element types.
Eventually, the capability of the method to accurately compute quantities of interest depending on the solution of the Stokes equation (e.g. the drag force) is discussed and Section~\ref{sc:Conclusion} summarizes the results of this paper.

\section{Stokes flow with strongly enforced symmetry of the stress tensor}
\label{sc:StokesSym}

In this section, the framework to handle symmetric tensors discussed by Fish and Belytschko in \cite{FishBelytschko2007} is recalled and the governing equations of a Stokes flow are formulated using Voigt notation.
First, the Cauchy formulation of the Stokes equation is recalled.

\subsection{Cauchy formulation of the Stokes flow}
\label{sc:Stokes}

Consider an open bounded domain $\Omega\subset\mathbb{R}^{\nsd}$ with boundary $\partial\Omega = \Gamma_D \cup \Gamma_N$, $\Gamma_D \cap \Gamma_N = \emptyset$ and $\nsd$ being the number of spatial dimensions. 
The strong form of the problem under analysis reads as follows:
\begin{equation} \label{eq:Stokes}
\left\{\begin{aligned}
- \grad\cdot \stress &= \bm{s}       &&\text{in $\Omega$,}\\
\grad\cdot\bu &= 0  &&\text{in $\Omega$,}\\
\stress &= -p \Insd + 2 \nu \defo{\bu}  &&\text{in $\Omega$,}\\
\bu &= \bu_D  &&\text{on $\Gamma_D$,}\\
\bn \cdot \stress &= \bm{t}  &&\text{on $\Gamma_N$,}\\                                          
\end{aligned}\right.
\end{equation}
where the pair $(\bu,p)$ represents the velocity and pressure fields and $\stress$ is the Cauchy stress tensor.
The terms $\bm{s}$, $\bu_D$ and $\bm{t}$ respectively are the volumetric source term, the Dirichlet boundary datum to impose the value of the velocity on $\Gamma_D$ and the traction applied on the Neumann boundary $\Gamma_N$.
The third equation, known as Stokes law, provides the relationship between the stress tensor and the velocity and pressure variables, through the viscosity coefficient $\nu > 0$, the $\nsd \times \nsd$ identity matrix $\Insd$ and the strain rate tensor $\defo{\bu}$, $\defo{} := \frac{1}{2}\left( \grad + \grad^T \right)$ being the symmetric part of the gradient.

It is well-known that the Cauchy and the velocity-pressure formulations of the Stokes equation are equivalent from the variational point of view.
Nevertheless, a major difference arises when considering the imposition of Neumann boundary conditions. 
On the one hand, natural boundary conditions for the Cauchy formulation enforce the value of the normal stress which represents a physical traction. 
On the other hand, the velocity-pressure formulation only accounts for the gradient of the velocity field instead of its symmetric part, leading to the imposition of the so-called \emph{pseudo-tractions}.
Hence, the physical interpretation is lost \cite{donea2003finite}.
Within this context, an artificial handling of Neumann boundary conditions is required to impose physically meaningful tractions. 
This represents a drawback when dealing with real-life and industrial applications in which the enforcement of physically relevant quantities is a major constraint to perform reliable numerical simulations and compare them with experimental data.

\subsection{Voigt notation for symmetric tensors}
\label{sc:Voigt}

The so-called Voigt notation relies on the idea of storing a second-order tensor in a vectorial format by appropriately rearranging its diagonal and off-diagonal components.
Consequently, the application of differential operators (e.g. symmetric gradient, divergence and $\curl$) and the geometrical projections (e.g. in the normal and tagential directions to a surface) may be expressed as matrix equations.
For this purpose, the rationale for the construction of differential operator and geometrical quantities using Voigt notation is recalled.

Consider the previously defined strain rate tensor $\defo{\bu}$.
Owing to its symmetry, only $\msd = \nsd(\nsd+1)/2$ components (i.e. three in 2D and six in 3D) have to be stored and the following column vector in $\RR^{\msd}$ is obtained:
\begin{equation} \label{eq:strainVoigt}
\strainV :=\begin{cases}
\bigl[e_{11} ,\; e_{22} ,\; e_{12} \bigr]^T
&\text{in 2D,} \\
\bigl[e_{11} ,\; e_{22} ,\; e_{33} ,\; e_{12} ,\; e_{13} ,\; e_{23} \bigr]^T
&\text{in 3D.} 
\end{cases}
\end{equation}
The components of the strain rate in Equation~\eqref{eq:strainVoigt} read as
\begin{equation} \label{eq:NormalShear}
e_{ij} := \frac{\partial u_i}{\partial x_j} + (1-\delta_{ij}) \frac{\partial u_j}{\partial x_i}, \quad \text{for } i,j = 1,\dotsc,\nsd ,
\end{equation}
where $\delta_{ij}$ is the classical Kronecker delta.
In order to retrieve the aforementioned strain rate tensor $\defo{\bu}$, the off-diagonal terms $e_{ij}, \ i \neq j$ have to be multiplied by a factor $1/2$, namely
\begin{equation} \label{eq:strainClassicVoigt}
\defo{\bu} :=\begin{cases}
\begin{bmatrix}
e_{11} & e_{12}/2 \\
e_{12}/2 & e_{22}
\end{bmatrix} 
&\text{in 2D,} \\
\begin{bmatrix}
e_{11} & e_{12}/2 & e_{13}/2 \\
e_{12}/2 & e_{22} & e_{23}/2 \\
e_{13}/2 & e_{23}/2 & e_{33}
\end{bmatrix}
&\text{in 3D.} 
\end{cases}
\end{equation}

Similarly, the symmetry of the stress tensor $\stress$ is exploited to store only $\msd$ components in the column vector
\begin{equation} \label{eq:stressVoigt}
\stressV :=\begin{cases}
\bigl[\sigma_{11} ,\; \sigma_{22} ,\; \sigma_{12} \bigr]^T
&\text{in 2D,} \\
\bigl[\sigma_{11} ,\; \sigma_{22} ,\; \sigma_{33} ,\; \sigma_{12} ,\; \sigma_{13} ,\; \sigma_{23} \bigr]^T
&\text{in 3D.} 
\end{cases}
\end{equation}

\subsubsection{Differential operators using Voigt notation}
\label{sc:VoigtDiff}

Following \cite{FishBelytschko2007}, the strain rate tensor can be written as $\strainV = \gradS \bu$ by introducing the $\msd \times \nsd$ matrix 
\begin{equation} \label{eq:symmGrad}
\gradS :=\begin{cases}
\begin{bmatrix}
\partial/\partial x_1 & 0 & \partial/\partial x_2 \\
0 & \partial/\partial x_2 & \partial/\partial x_1
\end{bmatrix}^T
&\text{in 2D,} \\
\begin{bmatrix}
\partial/\partial x_1 & 0 & 0 & \partial/\partial x_2 & \partial/\partial x_3 & 0 \\
0 & \partial/\partial x_2 & 0 & \partial/\partial x_1 & 0 & \partial/\partial x_3 \\
0 & 0 & \partial/\partial x_3 & 0 & \partial/\partial x_1 & \partial/\partial x_2
\end{bmatrix}^T
&\text{in 3D.} 
\end{cases}
\end{equation}

As previously done for the strain rate tensor by introducing the matrix $\gradS$ accounting for the symmetric part of the gradient, the vorticity vector is handled through its skew-symmetric part.
More precisely, the vorticity $\vort := \grad \times \bu$ may be expressed in terms of Voigt notation as $\vort = \gradW \bu$ through the $\nrr \times \nsd$ matrix $\gradW$, with $\nrr=\nsd(\nsd-1)/2$ being the number of rigid body rotations in the space (i.e. one in 2D and three in 3D):
\begin{equation} \label{eq:curlVoigt}
\gradW :=\begin{cases}
\bigl[-\partial/\partial x_2 ,\; \partial/\partial x_1 \bigr]
&\text{in 2D,} \\
\begin{bmatrix}
0 & -\partial/\partial x_3 & \partial/\partial x_2 \\
\partial/\partial x_3 & 0 & -\partial/\partial x_1 \\
-\partial/\partial x_2 & \partial/\partial x_1 & 0
\end{bmatrix}
&\text{in 3D.} 
\end{cases}
\end{equation}
\begin{remark}
The $\curl$ of a vector $\bv$ in two dimensions is a scalar quantity.
Nevertheless, it can also be computed by embedding $\bv$ in the three dimensional space $\RR^3$ and setting its third component equal to zero. 
Within this contect, $\grad \times \bv$ may be interpreted as a vector whose magnitude is given by the aforementioned matrix operation $\gradW \bv$ and pointing along the third direction.
\end{remark}

\subsection{Cauchy formulation of the Stokes flow using Voigt notation}
\label{sc:VoigtStokes}

Owing to the notation introduced in this section, the Stokes constitutive law may be expressed as $\stressV = -\bE p + \bD \gradS \bu$, where the vector $\bE \in \RR^{\msd}$ and the matrix $\bD \in \RR^{\msd \times \msd}$ read as
\begin{equation} \label{eq:EDVoigt}
\bE :=\begin{cases}
\bigl[1 ,\; 1 ,\; 0 \bigr]^T
&\text{in 2D,} \\
\bigl[1 ,\; 1 ,\; 1 ,\; 0 ,\; 0 ,\; 0 \bigr]^T
&\text{in 3D.} 
\end{cases}
\qquad
\bD :=\begin{cases}
\begin{bmatrix}
2\nu \Insd & \bm{0}_{\nsd \times 1} \\
\bm{0}_{\nsd \times 1}^T & \nu
\end{bmatrix}
&\text{in 2D,} \\
\begin{bmatrix}
2\nu \Insd & \bm{0}_{\nsd} \\
\bm{0}_{\nsd} & \nu \Insd
\end{bmatrix}
&\text{in 3D.} 
\end{cases}
\end{equation}

Moreover, the Neumann boundary condition applied on $\Gamma_N$ can be written as $\bN^T \stressV = \bt$ by introducing the $\msd \times \nsd$ matrix
\begin{equation} \label{eq:normalVoigt}
\bN :=\begin{cases}
\begin{bmatrix}
n_1 & 0 & n_2 \\
0 & n_2 & n_1
\end{bmatrix}^T
&\text{in 2D,} \\
\begin{bmatrix}
n_1 & 0 & 0 & n_2 & n_3 & 0\\
0 & n_2 & 0 & n_1 & 0 & n_3 \\
0 & 0 & n_3 & 0 & n_1 & n_2
\end{bmatrix}^T
&\text{in 3D.} 
\end{cases}
\end{equation}
accounting for the normal direction to the boundary.

Similarly, the projection of a vector along the tangential direction $\btau$, namely a tangent line in 2D and a tangent surface in 3D, reads as $\bu \cdot \btau = \bu^T \bT$, being $\bT \in \RR^{\nsd \times \nrr}$ the matrix
\begin{equation} \label{eq:tangentVoigt}
\bT :=\begin{cases}
\bigl[n_2 ,\; -n_1 \bigr]^T
&\text{in 2D,} \\
\begin{bmatrix}
0 & -n_3 & n_2 \\
n_3 & 0 & -n_1 \\
-n_2 & n_1 & 0
\end{bmatrix}
&\text{in 3D.} 
\end{cases}
\end{equation}

In order to rewrite Equation~\eqref{eq:Stokes} using Voigt notation, the divergence of a symmetric tensor is expressed in terms of the transpose of the matrix $\gradS$ accounting for the symmetric part of the gradient~\cite{FishBelytschko2007}.
In a similar fashion, recall that $\grad\cdot\bu = \tr(\grad\bu)$ and observe that the trace operator may be expressed via the vector $\bE$ introduced in Equation~\eqref{eq:EDVoigt}.
Combining the matrix forms of the symmetric gradient, the Stokes law and the normal direction presented above, the following formulation of the Stokes equation using Voigt notation is obtained:
\begin{equation} \label{eq:stokesSystemVoigt}
\left\{\begin{aligned}
-\gradS^T \stressV  &= \bm{s}       &&\text{in $\Omega$,}\\
\bE^T \gradS \bu &= 0        &&\text{in $\Omega$,}\\
\stressV &= -\bE p + \bD \gradS \bu        &&\text{in $\Omega$,}\\
\bu &= \bu_D  &&\text{on $\Gamma_D$,}\\
\bN^T \stressV &= \bt        &&\text{on $\Gamma_N$.}\\
\end{aligned}\right.
\end{equation}

\subsection{Fundamental theorems using Voigt notation}
\label{sc:VoigtTheo}

In \cite{preprintVoigtElasticity}, a generalized version of the Gauss's and Stokes' theorems using Voigt notation has been introduced. 
In order to construct the variational formulation of the problem under analysis, the following two lemmas are recalled.
\begin{lemma}[Generalized Gauss's theorem] \label{theo:Gauss}
	Consider a vector $\bv \in \RR^{\nsd}$ and a symmetric $\nsd \times \nsd$ tensor $\bTens$ whose counterpart in Voigt notation is $\bTensV$. It holds:
	\begin{equation} \label{eq:GenGauss}
	\int_{\partial\Omega} \left( \bN^T \bTensV \right) \cdot \bv \ d\Gamma = \int_{\Omega} \bTensV \cdot \left( \gradS \bv \right) d\Omega + \int_{\Omega} \left( \gradS^T \bTensV \right) \cdot \bv \ d\Omega .
	\end{equation}
\end{lemma}
\begin{lemma}[Generalized Stokes' theorem] \label{theo:Stokes}
	Consider a vector $\bv \in \RR^{\nsd}$. It holds:
	\begin{equation} \label{eq:GenStokes}
	\int_{\Omega} \gradW \bv \ d\Omega = \int_{\partial\Omega} \bv^T \bT \ d\Gamma .
	\end{equation}
\end{lemma}
The proofs follow straightforwardly by rewriting \eqref{eq:GenGauss}-\eqref{eq:GenStokes} using the corresponding continuous differential operators, see \cite{preprintVoigtElasticity}.

\section{A hybridizable discontinuous Galerkin method}
\label{sc:HDG}

HDG is a discontinuous Galerkin method with hybridization based on a mixed formulation.
First, it is defined the so-called broken computational domain by introducing a partition of the domain $\Omega$ in $\numel$ disjoint subdomains $\Omega_e$ with boundaries $\partial\Omega_e$.
The internal interface $\Gamma$ reads as
\begin{equation} \label{eq:skeleton}
\Gamma := \left[ \bigcup_{e=1}^{\numel} \partial\Omega_e \right] \setminus \partial\Omega ,
\end{equation}
whereas the mesh skeleton is given by the union of internal and Neumann boundary faces, namely $\Gamma \cup \Gamma_N$.

In what follows, the classical $\eltwo$ internal products for vector-valued functions in $\Omega_e \subset \Omega$ and $\partial\Omega_e \subset \Gamma \cup \partial\Omega$ are considered:
\begin{equation} \label{eq:innerScalar}
(\bu,\bw)_{\Omega_e} := \int_{\Omega_e} \bu \cdot \bw \ d\Omega  , \qquad \langle \hat{\bu}, \hat{\bw} \rangle_{\partial\Omega_e} := \sum_{\Gamma_i \subset \partial\Omega_e} \int_{\Gamma_i} \hat{\bu} \cdot \hat{\bw} \ d\Gamma .
\end{equation}
Moreover, owing to the piecewise discontinuous nature of the functions involved in the HDG formulation, the \emph{jump} operator $\jump{\cdot}$ is defined along each portion of the interface as the sum of the values from the element on the right and the left, $\Omega_e$ and $\Omega_l$ \cite{AdM-MFH:08}:
\begin{equation}
\jump{\odot} = \odot_e + \odot_l .
\end{equation}

The second-order problem in Equation~\eqref{eq:stokesSystemVoigt} may thus be written as a system of first-order equations as follows:
\begin{equation} \label{eq:StokesBrokenFirstOrder}
\left\{\begin{aligned}
\bL + \bDHalf \gradS \bu &= \bm{0}    &&\text{in $\Omega_e$, and for $e=1,\dotsc ,\numel$,}\\	
\gradS^T \bigl( \bDHalf \bL + \bE \, p \bigr) &= \bm{s}          &&\text{in $\Omega_e$, and for $e=1,\dotsc ,\numel$,}\\
\bE^T \gradS \bu &= 0         &&\text{in $\Omega_e$, and for $e=1,\dotsc ,\numel$,}\\
\bu &= \bu_D     &&\text{on $\Gamma_D$,}\\
\bN^T (\bDHalf \bL + \bE \, p) &= -\bt         &&\text{on $\Gamma_N$,}\\
\jump{\bu \otimes \bn} &= \bm{0}  &&\text{on $\Gamma$,}\\
\jump{\bN^T (\bDHalf \bL + \bE \,p)} &= \bm{0}  &&\text{on $\Gamma$,}\\
\end{aligned} \right.
\end{equation}
where $\bL$ is the so-called mixed variable and the last two equations are the \emph{transmission conditions} enforcing the continuity of respectively the velocity and the flux across the interface $\Gamma$. 
\begin{remark}
In the case of purely Dirichlet boundary conditions (i.e. $\Gamma_N = \emptyset$), an additional constraint is required to avoid the indeterminacy of the pressure. 
A common choice relies on imposing zero mean value of the pressure on the boundary (cf. e.g. \cite{MR2485446,MR3209969,Nguyen-NPC:10}):
\begin{equation} \label{eq:meanPressure}
\frac{1}{|\partial\Omega|} \langle p, 1\rangle_{\partial\Omega} = 0 .
\end{equation}
\end{remark}

\subsection{Strong form of the local and global problems}
\label{sc:HDGstrong}

In a series of papers by Cockburn and co-workers~\cite{MR2485446,Nguyen-NPC:10,Nguyen-CNP:10,MR2772094}, the hybridizable discontinuous Galerkin formulation for Stokes flow has been theoretically and numerically analyzed.
Starting from the mixed formulation on the broken computational domain in Equation~\eqref{eq:StokesBrokenFirstOrder}, HDG features two stages.

First, a set of $\numel$ local problems are defined element-by-element to compute $(\bL_e,\bu_e,p_e)$ for $e=1,\dotsc ,\numel$:
\begin{equation} \label{eq:StokesStrongLocal}
\left\{\begin{aligned}
\bL_e + \bDHalf \gradS \bu_e &= \bm{0}    &&\text{in $\Omega_e$}\\	
\gradS^T \bDHalf \bL_e + \gradS^T \bE \, p_e &= \bm{s}          &&\text{in $\Omega_e$}\\
\bE^T \gradS \bu_e &= 0           &&\text{in $\Omega_e$}\\
\bu_e &= \bu_D     &&\text{on $\partial\Omega_e \cap \Gamma_D$,}\\
\bu_e &= \bhu  &&\text{on $\partial\Omega_e \setminus \Gamma_D$,}\\
\end{aligned} \right.
\end{equation}
where $\bhu$ is an independent variable representing the trace of the velocity on the mesh skeleton $\Gamma \cup \Gamma_N$.
Remark that Equation~\eqref{eq:StokesStrongLocal} is a purely Dirichlet boundary value problem. As previously observed, an additional constraint has to be added to remove the indeterminacy of the pressure, namely
\begin{equation} \label{eq:constraintLoc}
\frac{1}{|\partial\Omega_e|} \langle p_e, 1 \rangle_{\partial\Omega_e} = \rho_e ,
\end{equation}
where $\rho_e$ denotes the mean pressure on the boundary of the element $\Omega_e$.
Hence, for $e=1,\dotsc ,\numel$ the local problem in Equation~\eqref{eq:StokesStrongLocal} provides $(\bL_e,\bu_e,p_e)$ in terms of the global unknowns $\bhu$ and $\rho$.

The trace of the velocity $\bhu$ and the mean pressure $\rho$ on the element boundaries are determined by solving the global problem accounting for the transmission conditions and the Neumann boundary condition:
\begin{equation} \label{eq:StokesStrongGlobal}
\left\{\begin{aligned}
\jump{\bu \otimes \bn} &= \bm{0}  &&\text{on $\Gamma$,}\\
\jump{\bN^T (\bDHalf \bL + \bE \,p)} &= \bm{0}  &&\text{on $\Gamma$,}\\
\bN^T (\bDHalf \bL + \bE \,p) &= -\bt         &&\text{on $\Gamma_N$.}\\
\end{aligned} \right.
\end{equation}
The first equation is automatically satisfied due to the Dirichlet boundary condition $\bu_e = \bhu$ imposed in the local problems and the unique definition of the hybrid variable $\bhu$ on each face of the mesh skeleton.
Moreover, the divergence-free condition in the local problem induces the following compatibility condition for each element $\Omega_e, \ e=1,\dotsc ,\numel$
\begin{equation}\label{eq:divergenceFreeConstraint}
\langle \bhu \cdot \bn_e , 1 \rangle_{\partial \Omega_e \setminus \Gamma_D} + \langle \bu_D \cdot \bn_e, 1 \rangle_{\partial \Omega_e\cap \Gamma_D} = 0 .
\end{equation}
Consider the Voigt counterpart $\bE^T \gradS \bu_e = 0 $ of the aforementioned constraint (cf. Equation~\eqref{eq:StokesStrongLocal}).
The resulting compatibility condition reads as
\begin{equation}\label{eq:divergenceFreeVoigt}
\langle \bE^T \bN_e \bhu , 1 \rangle_{\partial \Omega_e \setminus \Gamma_D} + \langle \bE^T \bN_e \bu_D , 1 \rangle_{\partial \Omega_e\cap \Gamma_D} = 0 \quad \text{for $e=1,\dotsc ,\numel$}
\end{equation}
and it is utilized to close the global problem.

\subsection{Weak form of the local and global problems}
\label{sc:HDGweak}

Consider the following discrete functional spaces according to the notation introduced in \cite{RS-AH:16}:
\begin{subequations}\label{eq:HDG-Spaces}
	\begin{align} 
	\Vh(\Omega) & := \left\{ v \in \eltwo(\Omega) : v \vert_{\Omega_e} \in \Pk(\Omega_e) \;\forall\Omega_e \, , \, e=1,\dotsc ,\numel \right\} , \label{eq:spaceScalarElem} \\
	\VhHat(S) & := \left\{ \hv \in \eltwo(S) : \hv\vert_{\Gamma_i}\in \Pk(\Gamma_i)	\;\forall\Gamma_i\subset S\subseteq\Gamma\cup\partial\Omega \right\}, \label{eq:spaceScalarFace}
	\end{align}
\end{subequations}
where $\mathcal{P}^{k}(\Omega_e)$ and $\mathcal{P}^{k}(\Gamma_i)$ are the spaces of polynomial functions of complete degree at most $k$ in $\Omega_e$ and on $\Gamma_i$, respectively. 

The discrete weak formulation of the local problems in Equation~\eqref{eq:StokesStrongLocal} is as follows: for $e=1,\dotsc ,\numel$, given $\bu_D$ on $\Gamma_D$ and $\bhu^h$ on $\Gamma\cup\Gamma_N$, find $(\bL^h_e ,\bu^h_e,p^h_e) \in [\Vh(\Omega_e)]^{\msd} \times [\Vh(\Omega_e)]^{\nsd} \times \Vh(\Omega_e)$ such that
\begin{subequations}\label{eq:HDGStokesWeakLocalPre}
	\begin{align}
	&
	\begin{aligned}
		- (\bv,\bL^h_e)_{\Omega_e} + (\gradS^T & \bDHalf \bv, \bu^h_e)_{\Omega_e} = \\
		& \langle \bN_e^T \bDHalf \bv , \bu_D\rangle_{\partial\Omega_e\cap\Gamma_D} + \langle \bN_e^T \bDHalf \bv , \bhu^h \rangle_{\partial\Omega_e\setminus\Gamma_D} , 
	\end{aligned} \label{eq:HDGStokesWeakLocalLPre}
	\\
	&
	\begin{aligned}
		-(\gradS \bw, \bDHalf \bL^h_e)_{\Omega_e}  - (\bE^T & \gradS \bw, p^h_e)_{\Omega_e} \\
		& + \langle \bw ,\bN_e^T \widehat{\bigl( \bDHalf \bL^h_e + \bE \, p^h_e \bigr)} \rangle_{\partial\Omega_e} =  (\bw,\bm{s})_{\Omega_e} , 
	\end{aligned} \label{eq:HDGStokesWeakLocalUPre}
	\\
		&
	(\gradS^T \bE \, q, \bu^h_e)_{\Omega_e}  = \langle q , \bE^T \bN_e \bu_D \rangle_{\partial\Omega_e\cap\Gamma_D} + \langle q , \bE^T \bN_e \bhu^h \rangle_{\partial\Omega_e\setminus\Gamma_D} , \label{eq:HDGStokesWeakLocalPPre}
	\\
	&
	\frac{1}{|\partial\Omega_e|} \langle p^h_e, 1 \rangle_{\partial\Omega_e} = \rho^h_e ,
	\label{eq:HDGStokesWeakLocalConstraintPre}
	\end{align}
\end{subequations}
for all $(\bv ,\bw,q) \in [\Vh(\Omega_e)]^{\msd} \times [\Vh(\Omega_e)]^{\nsd} \times \Vh(\Omega_e)$.
The trace of the numerical normal flux in Equation~\eqref{eq:HDGStokesWeakLocalUPre} is defined as follows
\begin{equation} \label{eq:traceStokes}
\bN_e^T \widehat{\bigl(\bDHalf \bL^h_e + \bE \, p^h_e \bigr)} := 
\begin{cases}
\bN_e^T \bigl( \bDHalf \bL^h_e + \bE \, p^h_e \bigr) + \btau (\bu^h_e - \bu_D) & \text{on $\partial\Omega_e\cap\Gamma_D$,} \\
\bN_e^T \bigl( \bDHalf \bL^h_e + \bE \, p^h_e \bigr) + \btau (\bu^h_e - \bhu^h) & \text{elsewhere,}  
\end{cases}
\end{equation}
where the stabilization parameter $\btau$ plays a crucial role in the stability, accuracy and convergence properties of the resulting HDG method \cite{Jay-CGL:09,Nguyen-NPC:09,Nguyen-NPC:09b}.
By plugging Equation~\eqref{eq:traceStokes} into Equation~\eqref{eq:HDGStokesWeakLocalUPre} and integrating by parts, the symmetric form of the discrete weak local problem is obtained: for $e=1,\dotsc ,\numel$, given $\bu_D$ on $\Gamma_D$ and $\bhu^h$ on $\Gamma\cup\Gamma_N$, find $(\bL^h_e ,\bu^h_e,p^h_e) \in [\Vh(\Omega_e)]^{\msd} \times [\Vh(\Omega_e)]^{\nsd} \times \Vh(\Omega_e)$ that satisfy
\begin{subequations}\label{eq:HDGStokesWeakLocal}
	\begin{align}
	&
	\begin{aligned}
		- (\bv,\bL^h_e)_{\Omega_e} + (\gradS^T & \bDHalf \bv, \bu^h_e)_{\Omega_e} = \\
		& \langle \bN_e^T \bDHalf \bv , \bu_D\rangle_{\partial\Omega_e\cap\Gamma_D} + \langle \bN_e^T \bDHalf \bv , \bhu^h \rangle_{\partial\Omega_e\setminus\Gamma_D} , 
	\end{aligned} \label{eq:HDGStokesWeakLocalL}
	\\
	&
	\begin{aligned}
	(\bw, \gradS^T \bDHalf \bL^h_e)_{\Omega_e} + & \langle \bw, \btau \bu^h_e \rangle_{\partial\Omega_e} + (\bw, \gradS^T \bE \, p^h_e)_{\Omega_e} = \\
	 & (\bw,\bm{s})_{\Omega_e} + \langle \bw , \btau \bu_D\rangle_{\partial\Omega_e\cap\Gamma_D} + \langle \bw , \btau \bhu^h \rangle_{\partial\Omega_e\setminus\Gamma_D},
	\end{aligned} 
	 \label{eq:HDGStokesWeakLocalU}
	\\
		&
	(\gradS^T \bE \, q, \bu^h_e)_{\Omega_e}  = \langle q , \bE^T \bN_e \bu_D \rangle_{\partial\Omega_e\cap\Gamma_D} + \langle q , \bE^T \bN_e \bhu^h \rangle_{\partial\Omega_e\setminus\Gamma_D} , \label{eq:HDGStokesWeakLocalP}
	\\
	&
	\frac{1}{|\partial\Omega_e|} \langle p^h_e, 1 \rangle_{\partial\Omega_e} = \rho^h_e ,
	\label{eq:HDGStokesWeakLocalConstraint}
	\end{align}
\end{subequations}
for all $(\bv ,\bw,q) \in [\Vh(\Omega_e)]^{\msd} \times [\Vh(\Omega_e)]^{\nsd} \times \Vh(\Omega_e)$.
\begin{remark}
From a practical point of view, the constraint on the mean value of the pressure on the boundary of the element introduced in Equation~\eqref{eq:HDGStokesWeakLocalConstraint} is handled by means of a Lagrange multiplier.
Thus, the matrix associated with the resulting local problem has a saddle point structure~\cite{RS-SH:18}.
\end{remark}

For the global problem, the discrete weak formulation equivalent to \eqref{eq:StokesStrongGlobal} is: find $\bhu^h\in[\VhHat(\Gamma\cup\Gamma_N)]^{\nsd}$ and $\rho^h \in \RR^{\numel}$ such that
\begin{subequations} \label{eq:HDGStokesWeakGlobal}
\begin{align}
&
\begin{aligned}
\sum_{e=1}^{\numel}\Bigl\{
\langle \bhw, \bN_e^T \bDHalf \bL^h_e \rangle_{\partial\Omega_e\setminus\Gamma_D}
& + \langle \bhw, \bE^T \bN_e p^h_e \rangle_{\partial\Omega_e\setminus\Gamma_D} 
+ \langle \bhw,\btau\, \bu^h_e \rangle_{\partial\Omega_e\setminus\Gamma_D} \Bigr. \\
& \Bigl. 
- \langle \bhw,\btau\,\bhu^h \rangle_{\partial\Omega_e\setminus\Gamma_D}\Bigr\}
= -\sum_{e=1}^{\numel} \langle \bhw, \bt \rangle_{\partial\Omega_e\cap\Gamma_N},
\end{aligned}
\label{eq:HDGStokesWeakGlobalUHat}
\\
&
\langle \bE^T \bN_e \bhu , 1 \rangle_{\partial \Omega_e \setminus \Gamma_D} = - \langle \bE^T \bN_e \bu_D , 1 \rangle_{\partial \Omega_e\cap \Gamma_D} = 0 \quad \text{for} \ e=1,\dotsc, \numel ,
\label{eq:HDGStokesWeakGlobalConstraint}
\end{align}
\end{subequations}
for all $\bhw\in[\VhHat(\Gamma\cup\Gamma_N)]^{\nsd}$.

\subsection{Local post-process of the velocity field}
\label{sc:StokesPostProcess}

As usual in HDG, an element-by-element post-process procedure is considered to construct an improved approximation of the velocity field.
Modifying the Brezzi-Douglas-Marini (BDM) projection operator~\cite{brezzi1991mixed}, in \cite{MR2772094,Nguyen-CNP:10}, a technique to retrieve an $\Hdiv$-conforming and exactly divergence-free velocity field was discussed.
In this section, a simpler approach inspired by the work of Stenberg \cite{MR1035181} and exploited in \cite{Nguyen-NPC:10,RS-AH:16,RS-SH:18} is considered.
The requirement of $\Hdiv$-conformity is relaxed and the resulting local post-process problem exploits the optimal convergence rate of order $k+1$ of the mixed variable to construct a velocity field $\bu^\star$ superconverging with order $k+2$.

Nevertheless, it is known~\cite{Nguyen-CNP:10} that using the Cauchy formulation of the Stokes equation, a loss of superconvergence is experienced by low-order approximations.
The Voigt notation introduced in Section~\ref{sc:Voigt} allows to remedy this issue and to circumvent the complex mathematical framework of $\bm{M}$-decomposition discussed in \cite{MR3626530,MR3601005,MR3601012,MR3649424} to devise superconvergent HDG approximations with strongly and weakly symmetric stress tensors.
Following \cite{preprintVoigtElasticity}, the space $\Vh_\star(\Omega)$ of the polynomials of complete degree at most $k+1$ on each element $\Omega_e$
\begin{equation} \label{eq:Vstar}
\Vh_\star(\Omega) := \left\{ v \in \eltwo(\Omega) : v \vert_{\Omega_e}\in \mathcal{P}^{k+1}(\Omega_e) \;\forall\Omega_e \, , \, e=1,\dotsc ,\numel \right\} 
\end{equation}
is introduced.
For each element $\Omega_e, \ e=1,\dotsc,\numel$, the post-processed velocity $\bu^\star$ is the solution of the problem
\begin{equation} \label{eq:postVoigt}
\left\{\begin{aligned}
\gradS^T \bDHalf \gradS \bu_e^\star  &= - \gradS^T \bL^h_e       &&\text{in $\Omega_e$,}\\	
\bN_e^T \bDHalf \gradS \bu_e^\star &= - \bN_e^T \bL^h_e        &&\text{on $\partial\Omega_e$,}\\
\end{aligned}\right.
\end{equation}
in the space $\left[ \Vh_\star(\Omega) \right]^{\nsd}$.
The element-by-element problem in Equation~\eqref{eq:postVoigt} is obtained by the definition of the mixed variable in Equation~\eqref{eq:StokesStrongLocal} and exploits the naturally equilibrated fluxes as condition on the boundary of the element.

The solution of Equation~\eqref{eq:postVoigt} is determined up to rigid motions, namely $\nsd$ translations and $\nrr$ rotations, being $\nsd = 2$ and $\nrr = 1$ in 2D and $\nsd = \nrr = 3$ in 3D.
According to \cite{preprintVoigtElasticity}, a set of $\nsd + \nrr$ constraints is introduced to retrieve the uniqueness of the solution.
On the one hand, the indeterminacy due to the $\nsd$ rigid translational modes is resolved introducing the following constraint on the mean value of the velocity:
\begin{equation} \label{eq:postprocessMean}
( \bu^\star_e , 1 )_{\Omega_e} = ( \bu^h_e , 1)_{\Omega_e} .
\end{equation} 
On the other hand, the $\nrr$ rigid rotational modes are taken care of by means of a condition on the $\curl$ of the velocity, namely
\begin{equation} \label{eq:postprocessRot}
( \grad \times \bu^\star_e , 1 )_{\Omega_e}  =  \langle \bhu^h \cdot \btauE , 1 \rangle_{\partial \Omega_e} ,
\end{equation}
where the right-hand side of Equation~\eqref{eq:postprocessRot} follows from the application of Stokes' theorem, being $\bu^h_e = \bhu^h$ on $\partial\Omega_e$ and $\btauE$ the tangential direction to the boundary $\partial\Omega_e$.

It is worth noting that other conditions may be considered to resolve the indeterminacy of the problem in Equation~\eqref{eq:postVoigt}.
Nevertheless, in order for the post-processed velocity to be superconvergent, the quantities appearing in these constraints have to converge with order $m \geq k+2$.
If this is not the case, despite the resulting system admits a unique solution, the superconvergence property is lost.
For the strategy discussed in the present work, extensive numerical experiments have shown that the right-hand sides of both \eqref{eq:postprocessMean} and \eqref{eq:postprocessRot} converge with order $m > k+2$.
A rigorous proof of this result is currently under investigation.
\begin{remark}
Recall that the $\curl$ of the velocity represents the vorticity of the fluid.
Within this context, the left hand side of Equation~\eqref{eq:postprocessRot} may be physically interpreted as the mean value of the vorticity inside the element $\Omega_e$.
Similarly, the right-hand side represents the circulation of the flow around the boundary $\partial\Omega_e$.
\end{remark}
Eventually, by exploiting the Voigt notation, Equation~\eqref{eq:postprocessRot} is equivalent to
\begin{equation} \label{eq:postprocessVoigt}
( \gradW \bu^\star_e , 1 )_{\Omega_e} = \langle \bigl[ \bhu^h \bigr]^T \bT , 1 \rangle_{\partial \Omega_e} .
\end{equation}

\section{Numerical studies}
\label{sc:simulations}

In this section, several examples with known analytical solution are considered, in two and three dimensions, to verify the optimal convergence and superconvergence properties of the error of the primal, mixed and post-processed variables, measured in the $\eltwo(\Omega)$ norm and for different element types.
First, a numerical study of the influence of the stabilization parameter $\btau$ on the accuracy of the proposed HDG method is performed.

\subsection{Influence of the stabilization parameter}
\label{sc:stabilization}

As previously stated and extensively studied in a series of publications by Cockburn and co-workers (cf. e.g. \cite{Jay-CGL:09,Nguyen-NPC:09,Nguyen-NPC:09b}), the HDG stabilization parameter has an important effect on the convergence properties of the method.
For the sake of simplicity, a stabilization tensor of the form $\btau = \tau \Insd$, equal on all the faces of the internal skeleton $\Gamma \cup \Gamma_N$ is considered.
In what follows, a numerical study of the role of the scalar parameter $\tau$ is presented.

\subsubsection{Two dimensional example}
\label{sc:stabilization2D}

The first example considers the well-known problem of the Wang flow in the domain $\Omega = [0,1]^2$. 
The source term $\bm{s}$ is selected so that the analytical velocity field has the following expression
\begin{equation}
  \bu(\bx) =
  \begin{Bmatrix} 
  2 a x_2 - b \lambda \cos(\lambda x_1) \exp\{-\lambda x_2\} \\[1ex] 
  b \lambda \sin(\lambda x_1) \exp\{-\lambda x_2\}
  \end{Bmatrix} ,
\end{equation}
whereas the pressure is uniformly zero in the domain.
The values $a = b = \lambda = 1$ are set for the constants and the kinematic viscosity $\nu$ is taken equal to 1.
Neumann boundary conditions, corresponding to the analytical normal flux, are imposed on $\Gamma_N = \{(x_1,x_2) \in \Omega \; | \; x_2=0\}$ and the analytical velocity field is enforced on $\Gamma_D = \partial \Omega \setminus \Gamma_N$ via Dirichlet boundary conditions.

\begin{figure}[!tb]
	\centering
	\subfigure[Quadrilateral mesh]{\includegraphics[width=0.28\textwidth]{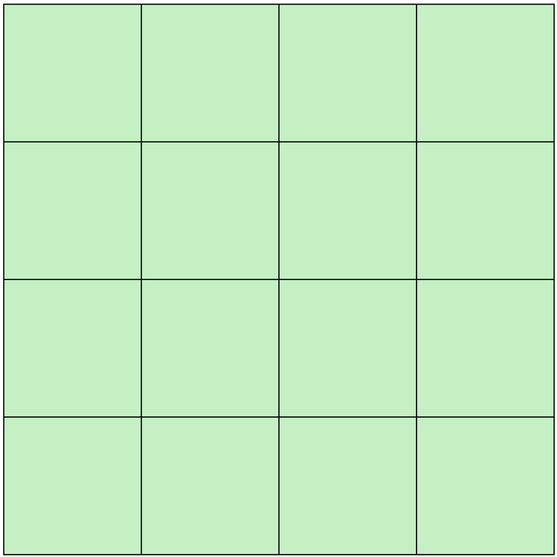}}
	\subfigure[Triangular mesh \#1]{\includegraphics[width=0.28\textwidth]{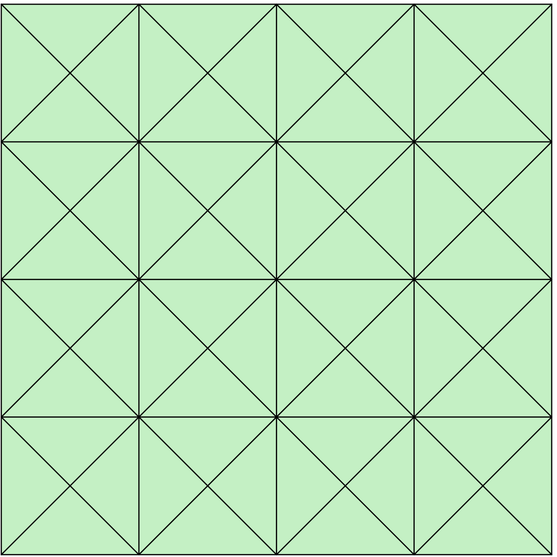}}
	\subfigure[Triangular mesh \#2]{\includegraphics[width=0.28\textwidth]{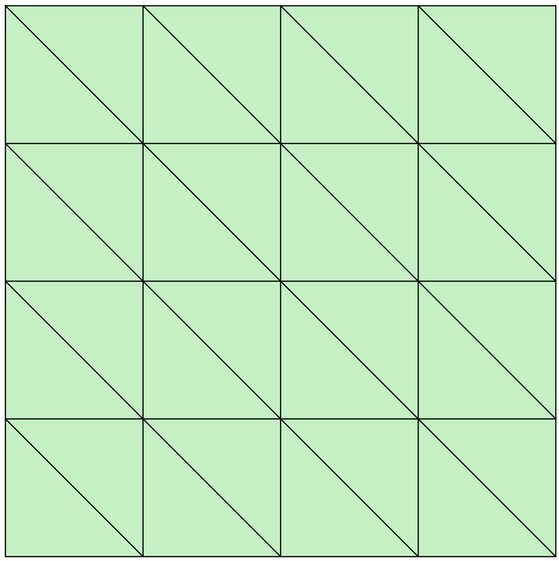}}
	
	\caption{Second level of refinement for three types of two dimensional meshes of $\Omega=[0,1]^2$ utilized for the mesh convergence study.}
	\label{fig:2Dmeshes}
\end{figure}

Uniform meshes of quadrilateral and triangular elements are considered. The second level of refinement of the meshes is shown in Figure~\ref{fig:2Dmeshes}.
It is worth noting that the triangular mesh \#1 has considerably more degrees of freedom than the triangular mesh \#2 for a similar characteristic size. 

\begin{figure}[!tb]
	\centering
	\subfigure[$u_1$]{\includegraphics[height=0.24\textwidth]{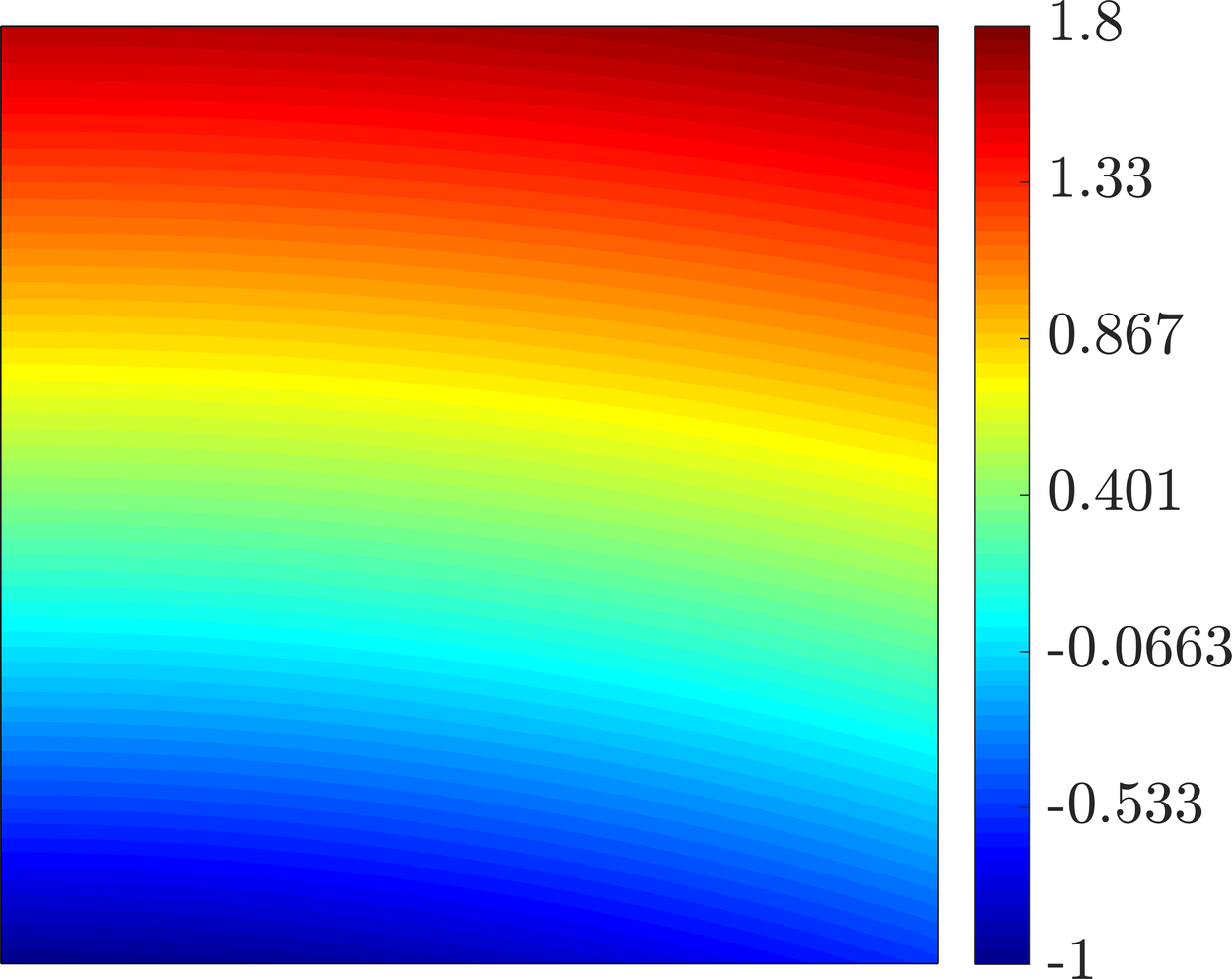}}
	\subfigure[$u_2$]{\includegraphics[height=0.24\textwidth]{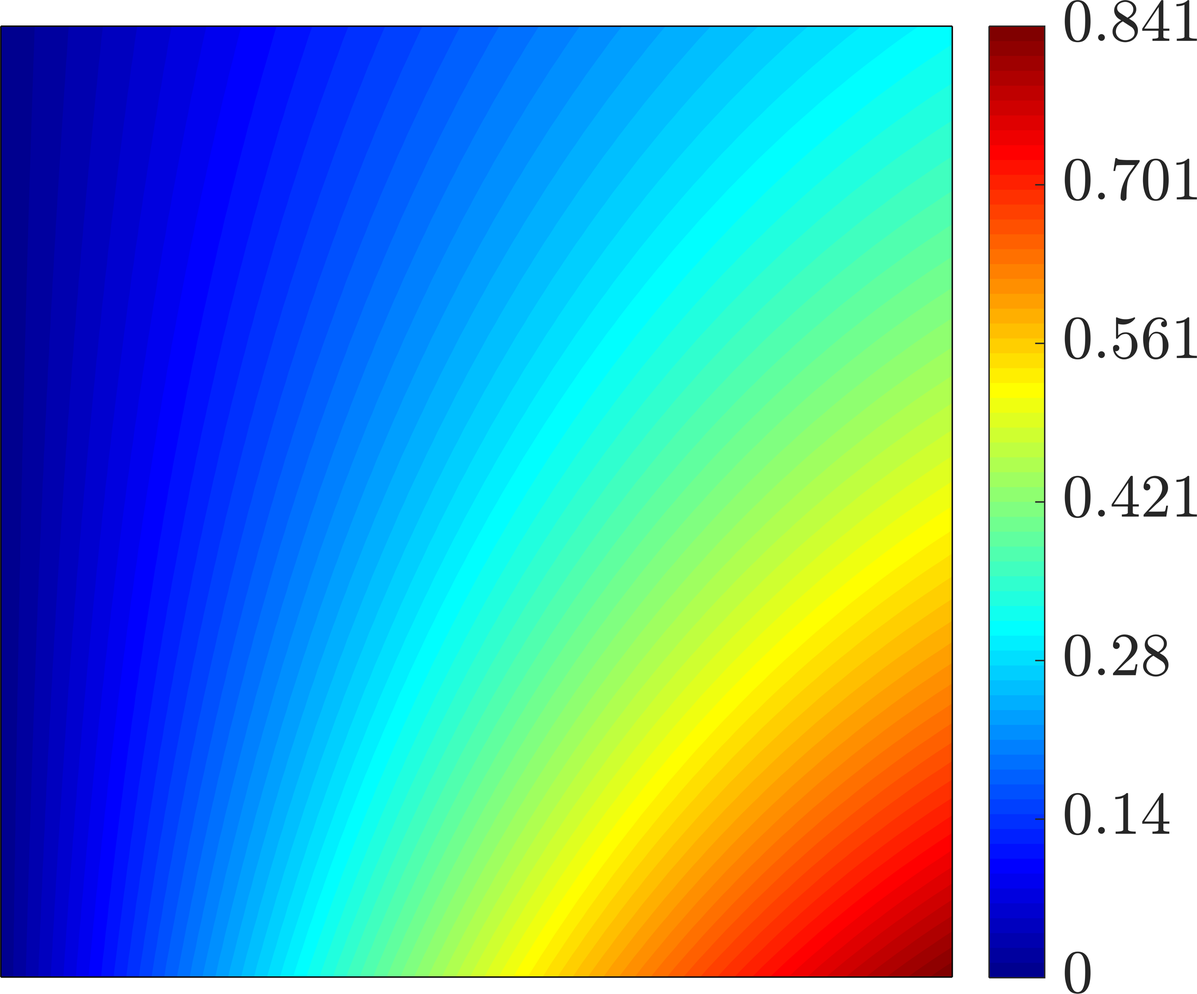}}
	
	\caption{Two dimensional problem: HDG approximation of the velocity field using the fourth refinement of the triangular mesh \#2 and $k=2$.}
	\label{fig:2Dsol}
\end{figure}

The components of the velocity field computed on the fourth level of refinement of the triangular mesh \#2 and using a quadratic degree of approximation are depicted in Figure~\ref{fig:2Dsol}.

\begin{figure}[!tb]
	\centering
	\subfigure[Quadrilateral mesh]{\includegraphics[width=0.4\textwidth]{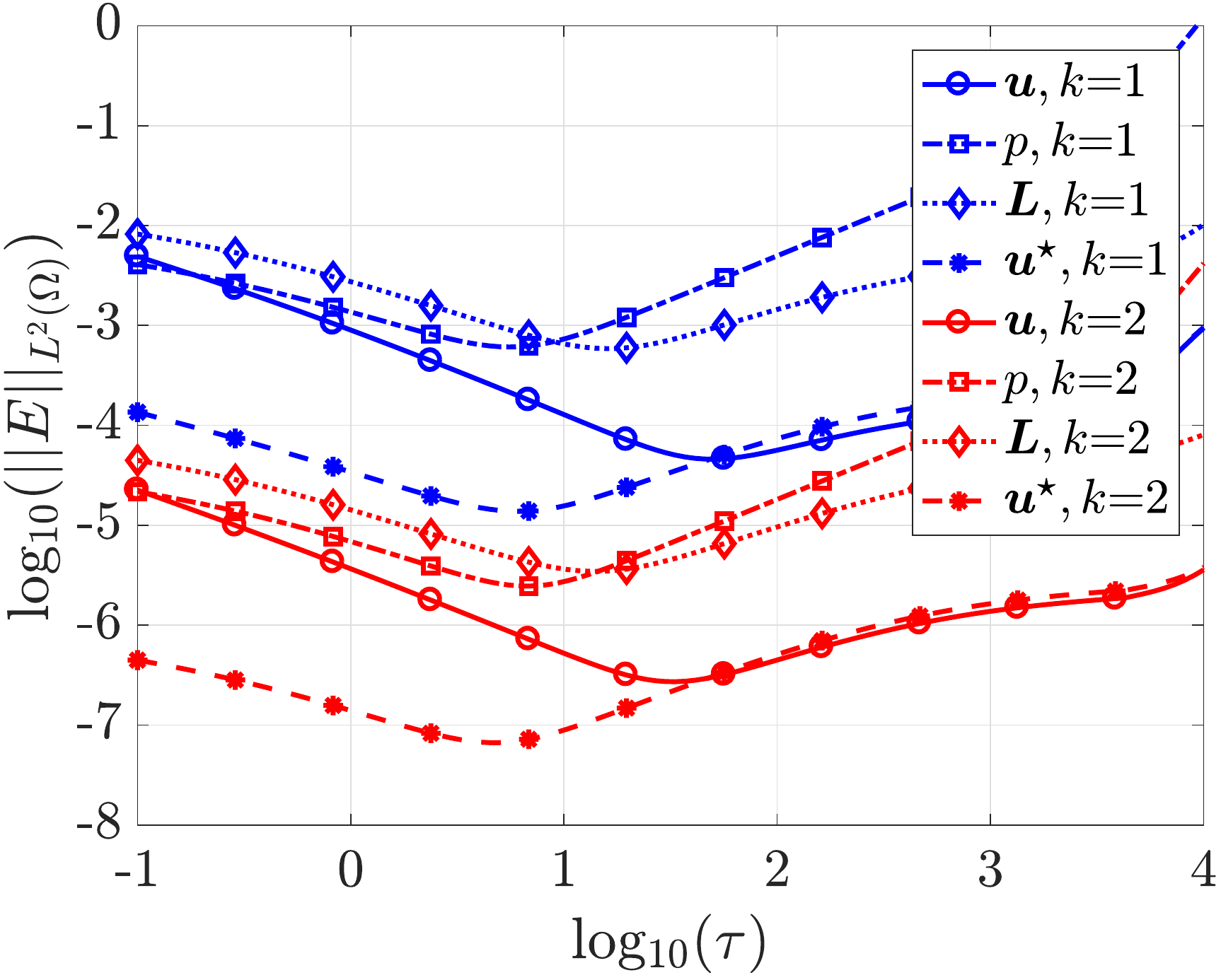}}
	\subfigure[Triangular mesh \#1]{\includegraphics[width=0.4\textwidth]{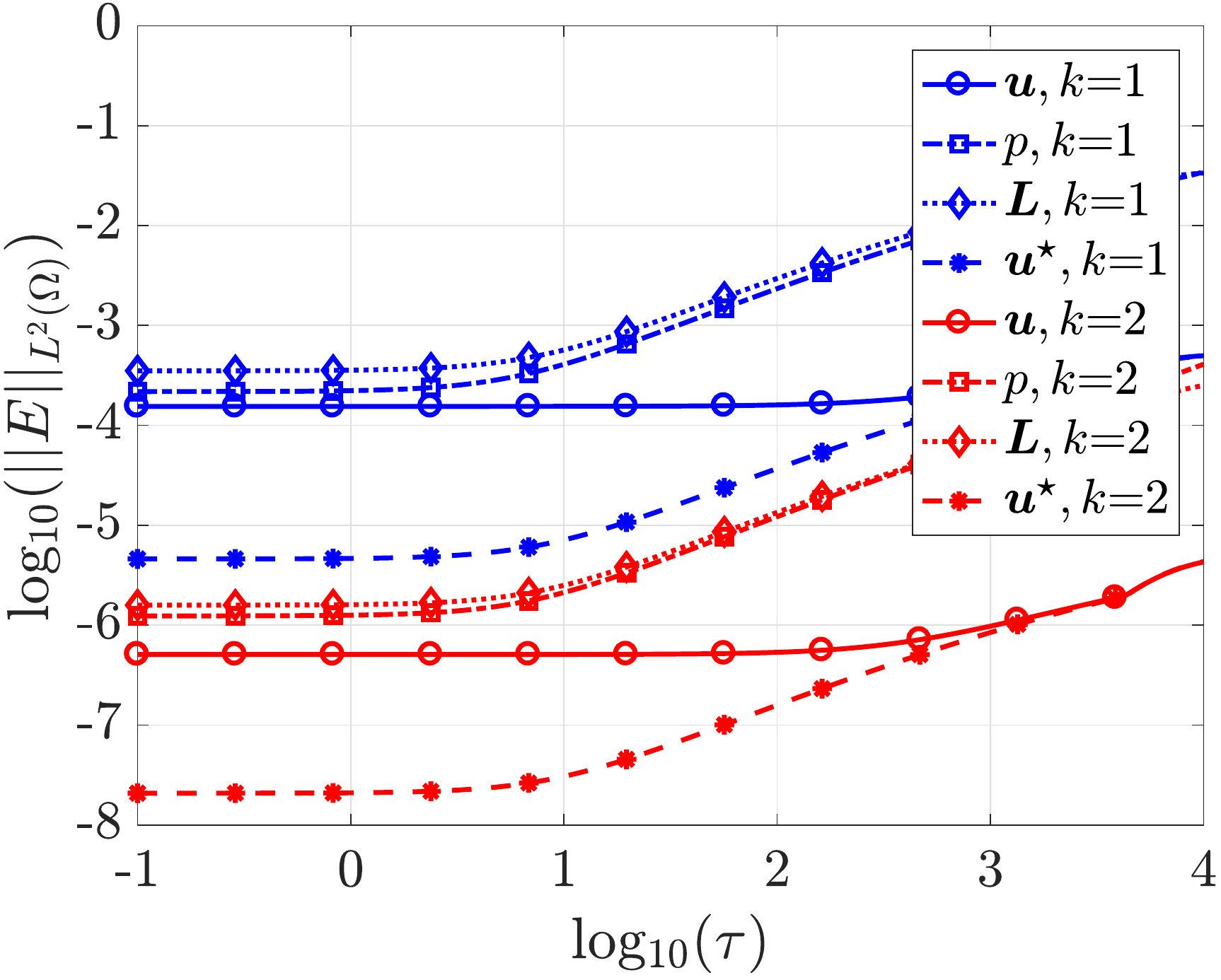}}	
	\subfigure[Triangular mesh \#2]{\includegraphics[width=0.4\textwidth]{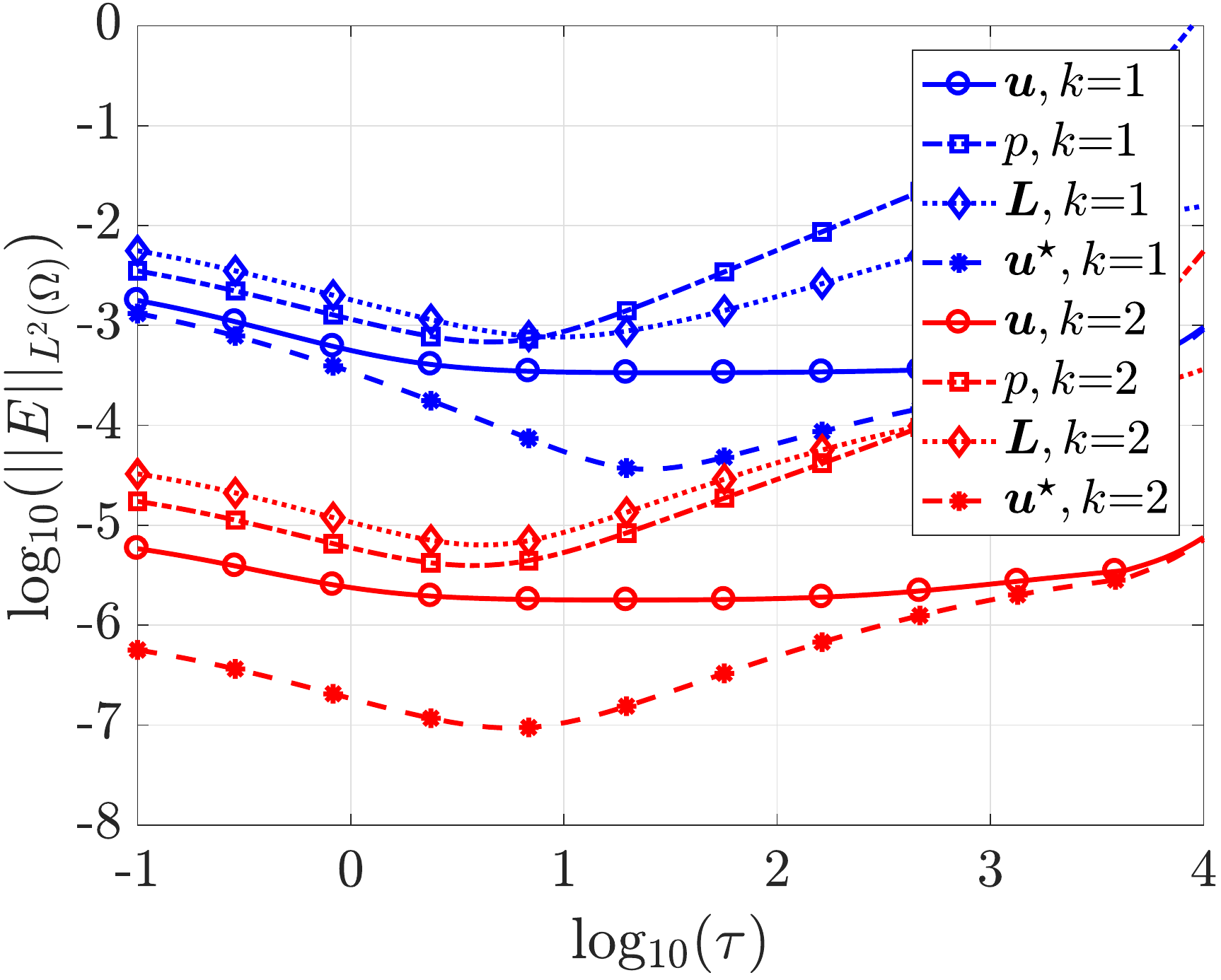}}	
	
	\caption{Two dimensional problem: error of the primal, mixed and post-processed variables, $\bu$, $p$, $\bL$ and $\bu^\star$, in the $\eltwo(\Omega)$ norm as a function of the stabilization parameter and for the fourth level of mesh refinement.}
	\label{fig:tauInfluence2D}
\end{figure}

Figure~\ref{fig:tauInfluence2D} shows the evolution of the error of the primal, mixed and post-processed variables, $\bu$, $p$, $\bL$ and $\bu^\star$, in the $\eltwo(\Omega)$ norm as a function of the stabilization parameter $\tau$.
The numerical study is performed on the fourth level of mesh refinement, using polynomial approximations of complete degree 1 and 2 and values of $\tau$ spanning from 0.1 to 10,000.
It is straightforward to observe that for all the meshes under analysis, there exists a value of $\tau$ minimizing the $\eltwo(\Omega)$ norm of the error of the velocity.
Nevertheless, to guarantee the accuracy of the approximation, the $\hone(\Omega)$ norm of the error should be accounted for and consequently both $\bu$ and $\bL$ are considered in the choice of the optimal value of $\tau$.
Within this context and in order for the post-processed velocity field $\bu^\star$ to provide a gain in accuracy with respect to $\bu$, the value $\tau = 4$ is chosen for quadrilateral meshes and triangular meshes of the first type.
For triangular meshes of the second type, the minimum of the error in the primal variable is achieved for values of $\tau$ substantially larger than 10.
Despite the approximation of the mixed variable deteriorates when the stabilization parameter increases, this effect is limited for values of $\tau < 50$.
The value of $\tau=40$ is thus considered as it provides a good compromise for the quality of the approximation of the primal, mixed and post-processed variables.
\begin{remark}
Consider the family of meshes in Figure~\ref{fig:2Dmeshes}.
The triangular mesh \#1 features one node located in the barycenter of each underlying quadrilateral. The resulting mesh provides significantly more information than the triangular mesh \#2 of the corresponding refinement level.
Thus, owing to the aforementioned extra node and to the tensorial nature of the basis functions defined on the quadrilateral meshes, the behavior of the triangular meshes \#1 is expected to be more similar to the quadrilateral ones than to the triangular meshes \#2, as observed in the previous numerical simulations in Figure~\ref{fig:tauInfluence2D}.
\end{remark}

\subsubsection{Three dimensional example}
\label{sc:stabilization3D}

The second example, inspired by \cite{FLD:FLD1650190502}, is an analytical solution of the problem in Equation~\eqref{eq:Stokes} set in the domain $\Omega = [0,1]^3$. 
The source term is selected so that the analytical velocity is 
\begin{equation}
  \bu(\bx) =
  \begin{Bmatrix} 
  b \exp\{a(x_1{-}x_3) + b(x_2{-}x_3)\} - a \exp\{ a(x_3{-}x_2) + b(x_1{-}x_2)\} \\[1ex] 
  b \exp\{a(x_2{-}x_1) + b(x_3{-}x_1)\} - a \exp\{a(x_1{-}x_3) + b(x_2{-}x_3)\} \\[1ex]
  b \exp\{a(x_3{-}x_2) + b(x_1{-}x_2)\} - a \exp\{a(x_2{-}x_1) + b(x_3{-}x_1)\} 
  \end{Bmatrix} 
\end{equation}
and the corresponding pressure field is
\begin{equation}
  p(\bx) = x_1 (1-x_1) .
\end{equation}
The values $a = 1$ and $b = 0.5$ are considered and the kinematic viscosity $\nu$ is taken equal to 1.
Neumann boundary conditions, corresponding to the analytical flux, are imposed on $\Gamma_N = \{(x_1,x_2,x_3) \in \Omega \; | \; x_3=0\}$ and the analytical velocity field is enforced on $\Gamma_D = \partial \Omega \setminus \Gamma_N$ via Dirichlet boundary conditions.

\begin{figure}[!tb]
	\centering
	\subfigure[Hexahedral mesh]{\includegraphics[width=0.24\textwidth]{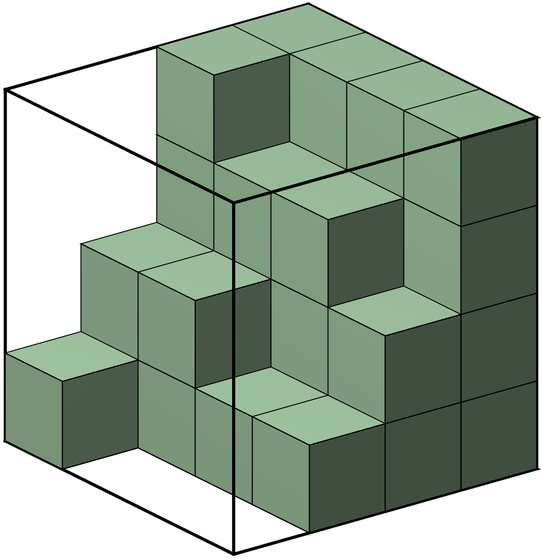}}
	\subfigure[Tetrahedral mesh]{\includegraphics[width=0.24\textwidth]{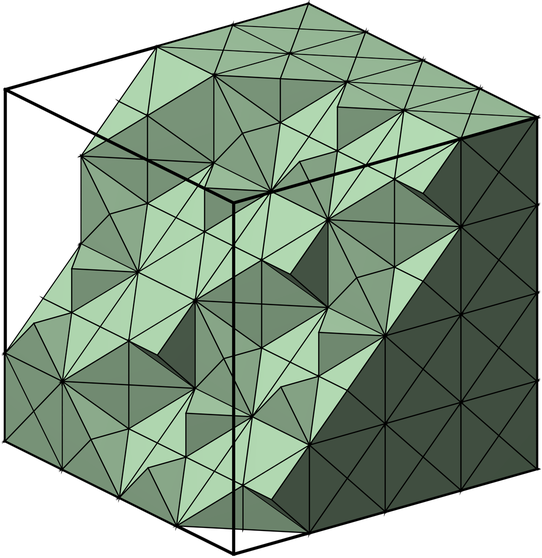}}		
	\subfigure[Prismatic mesh]{\includegraphics[width=0.24\textwidth]{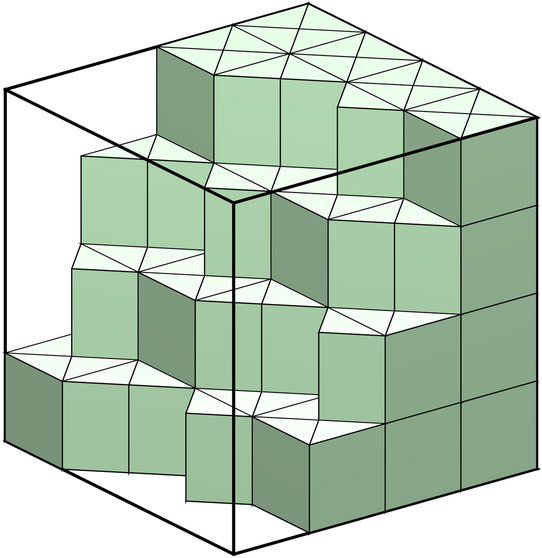}}
	\subfigure[Pyramidal mesh]{\includegraphics[width=0.24\textwidth]{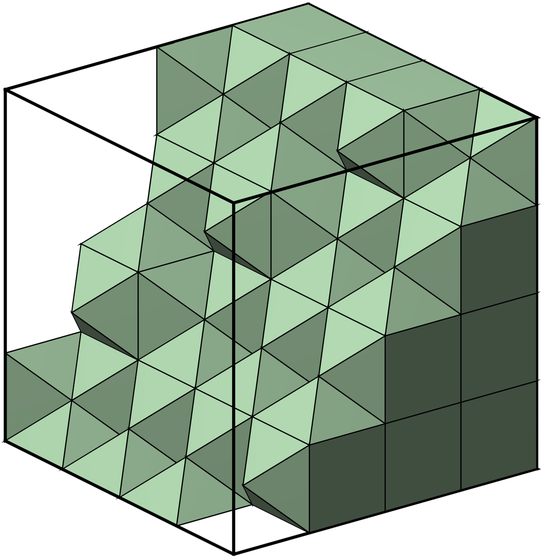}}
	
	\caption{Third level of refinement for four types of three dimensional meshes of $\Omega=[0,1]^3$ utilized for the mesh convergence study.}
	\label{fig:3Dmeshes}
\end{figure}

Figure~\ref{fig:3Dmeshes} shows a cut through the third level of refinement of the uniform meshes of hexahedral, tetrahedral, prismatic and pyramidal elements considered in this study.

\begin{figure}[!tb]
	\centering
	\subfigure[$u_1$]{\includegraphics[width=0.24\textwidth]{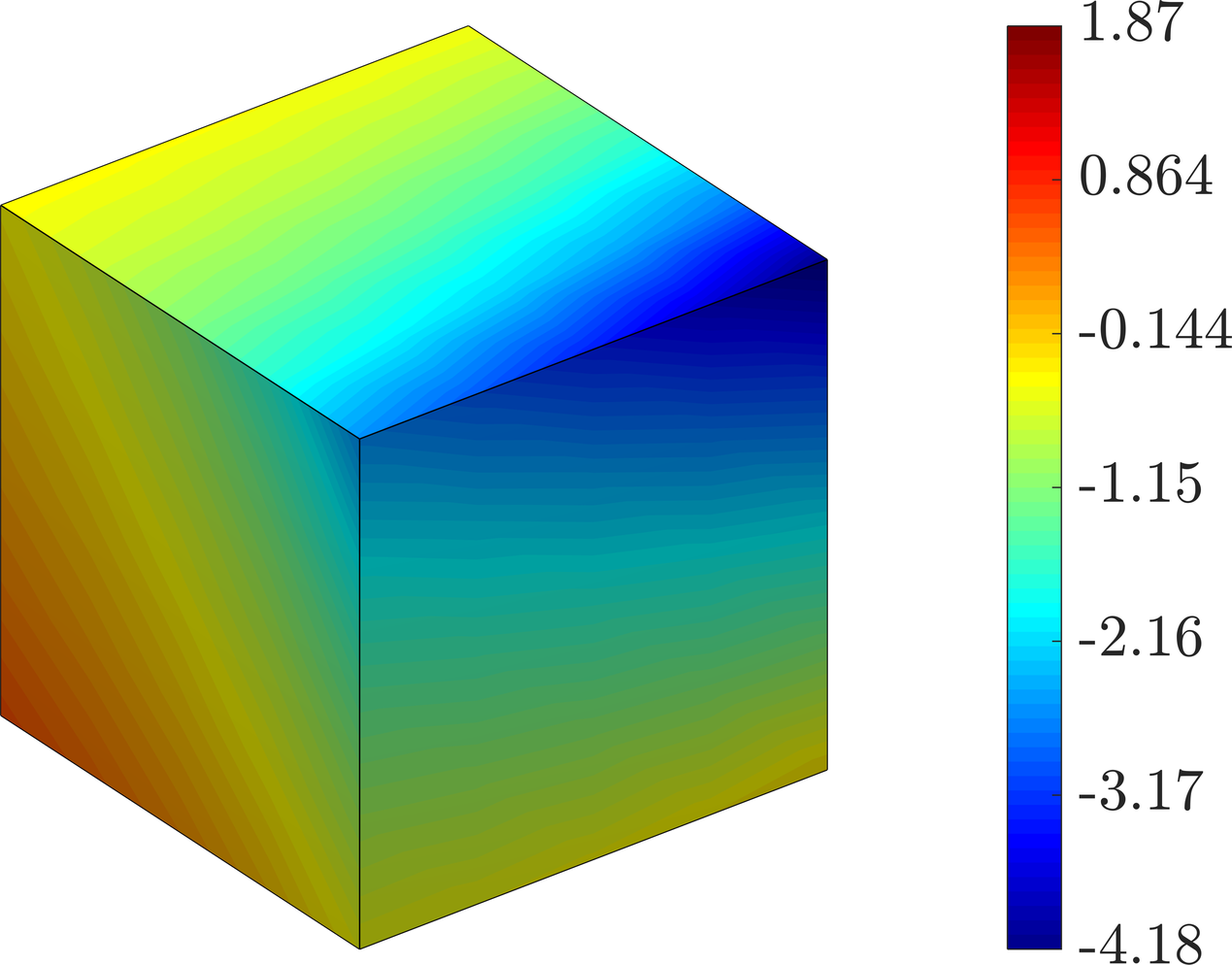}}
	\subfigure[$u_2$]{\includegraphics[width=0.24\textwidth]{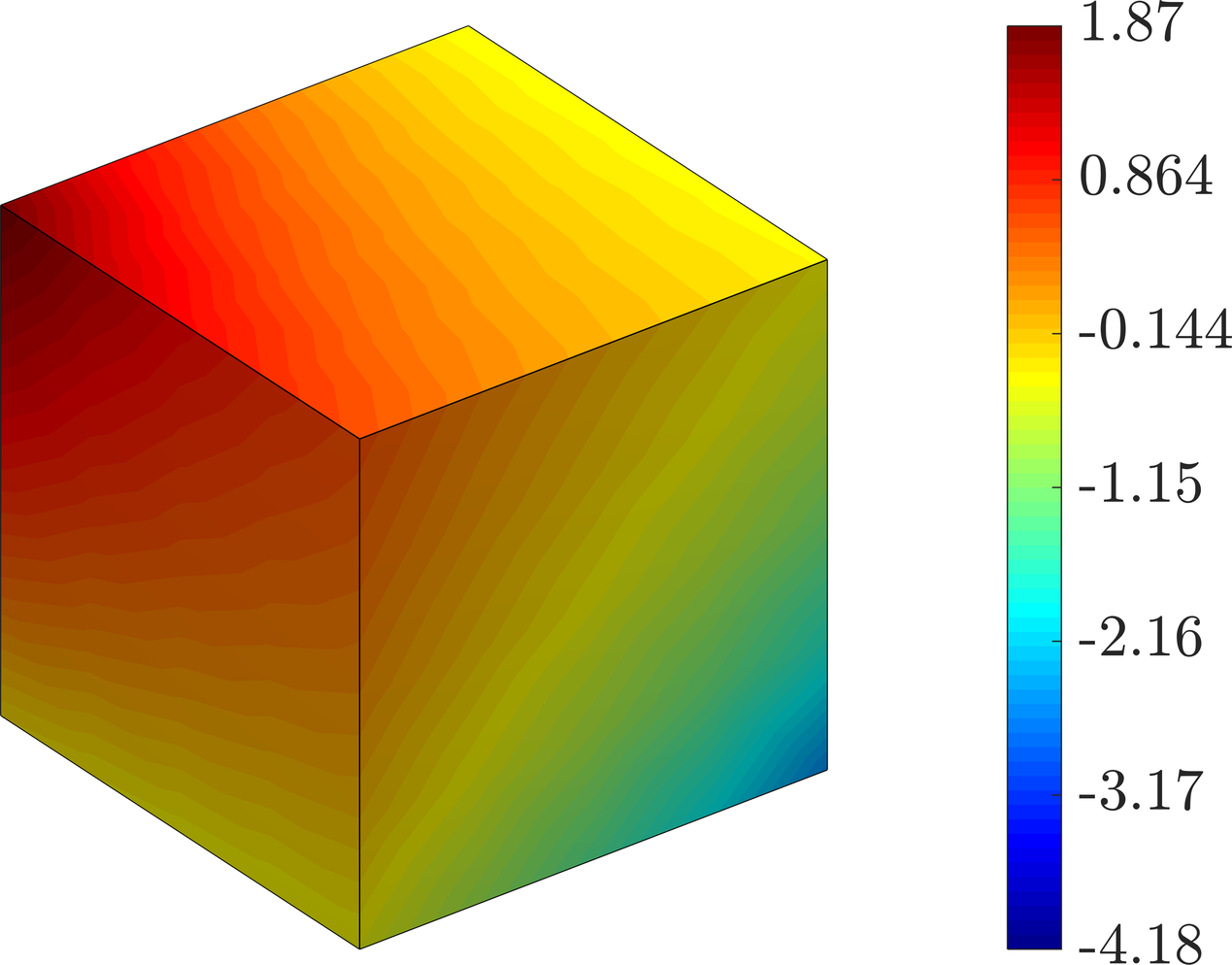}}
	\subfigure[$u_3$]{\includegraphics[width=0.24\textwidth]{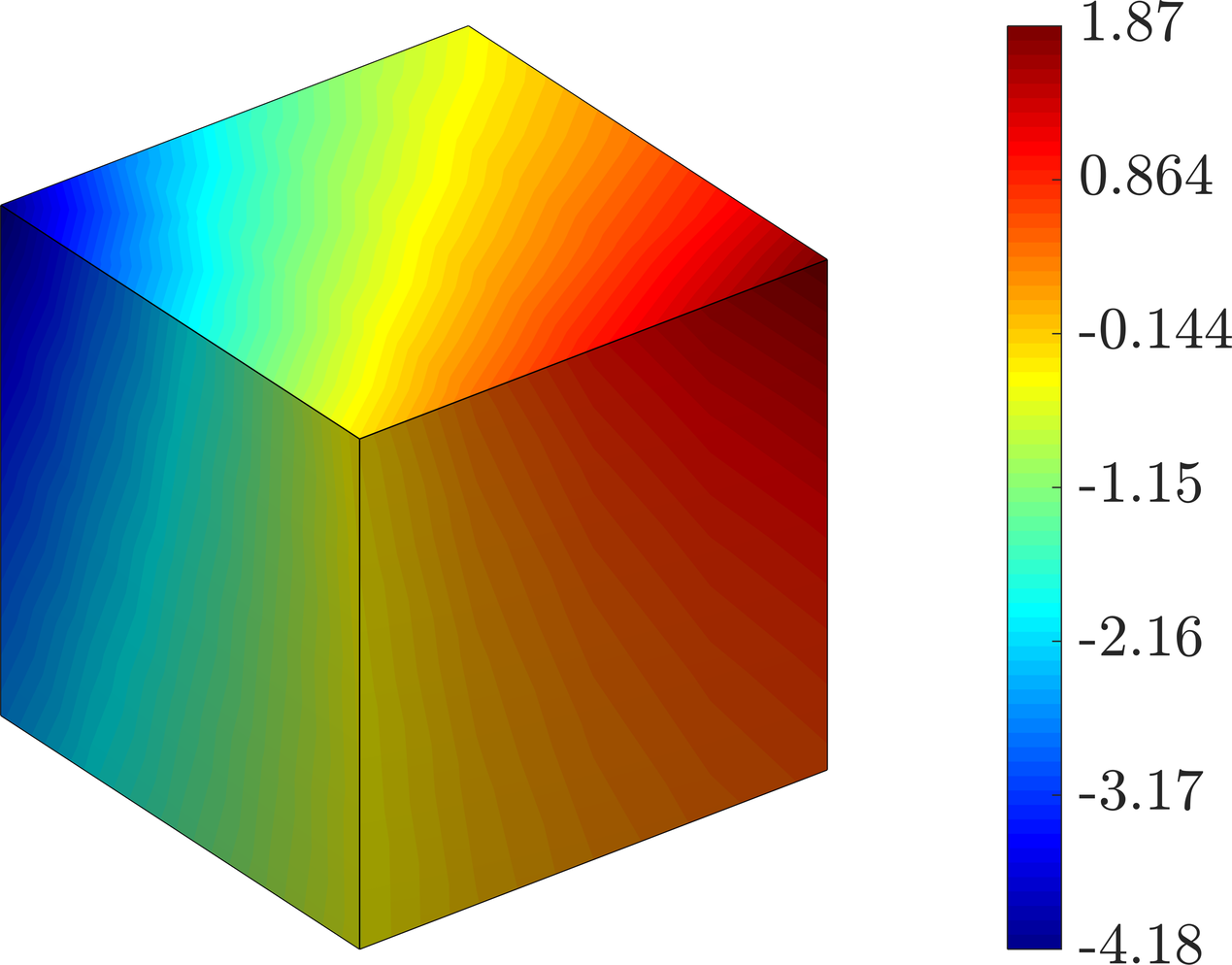}}
	\subfigure[$p$]{\includegraphics[width=0.24\textwidth]{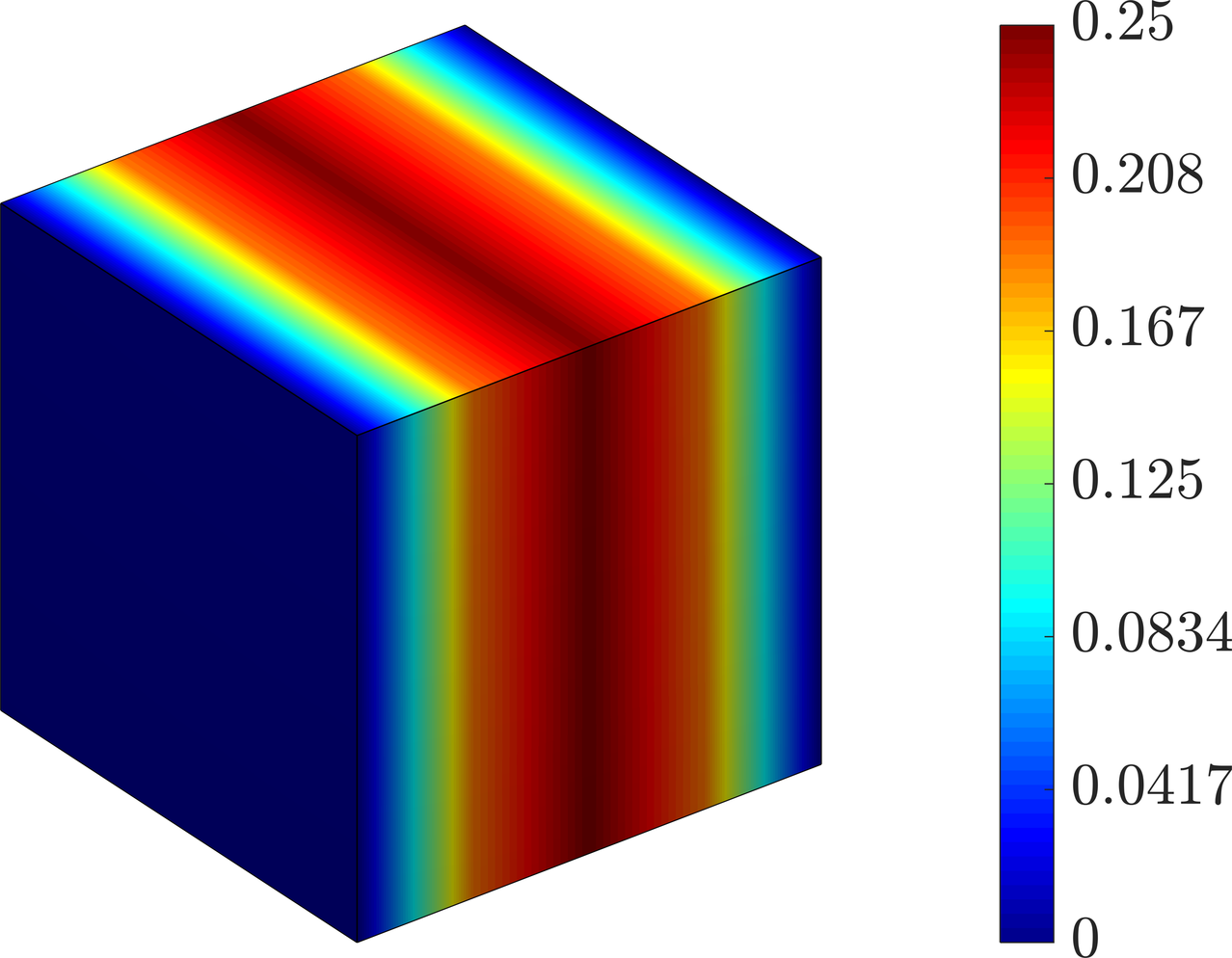}}
	
	\caption{Three dimensional problem: HDG approximation of the velocity and pressure fields using the third refinement of the hexahedral mesh and $k=3$.}
	\label{fig:3Dsol}
\end{figure}

The velocity and pressure fields computed on the third level of refinement of the hexahedral mesh and using a cubic degree of approximation are depicted in Figure~\ref{fig:3Dsol}.

\begin{figure}[!tb]
	\centering
	\subfigure[Hexahedral mesh]{\includegraphics[width=0.4\textwidth]{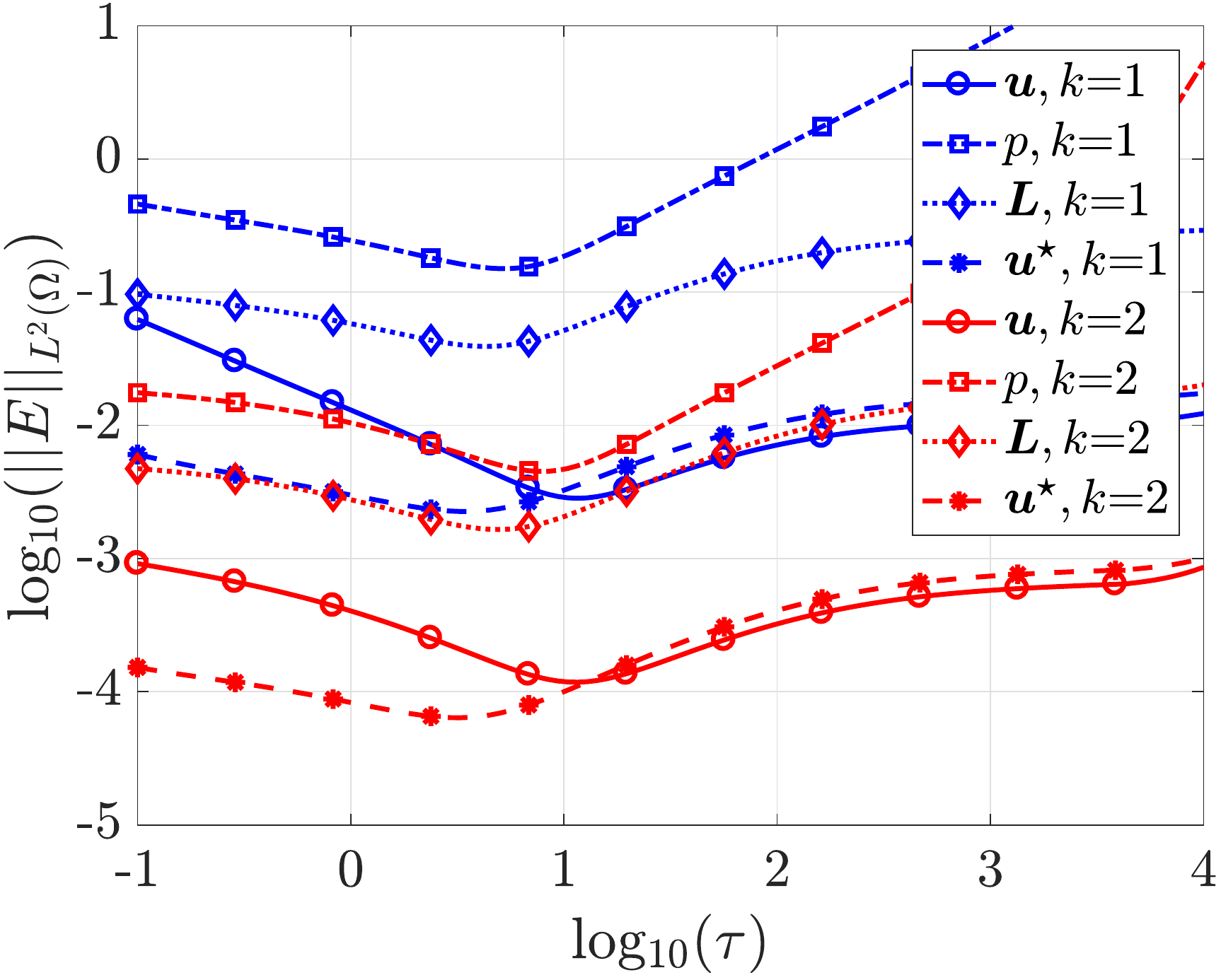}}
	\subfigure[Tetrahedral mesh]{\includegraphics[width=0.4\textwidth]{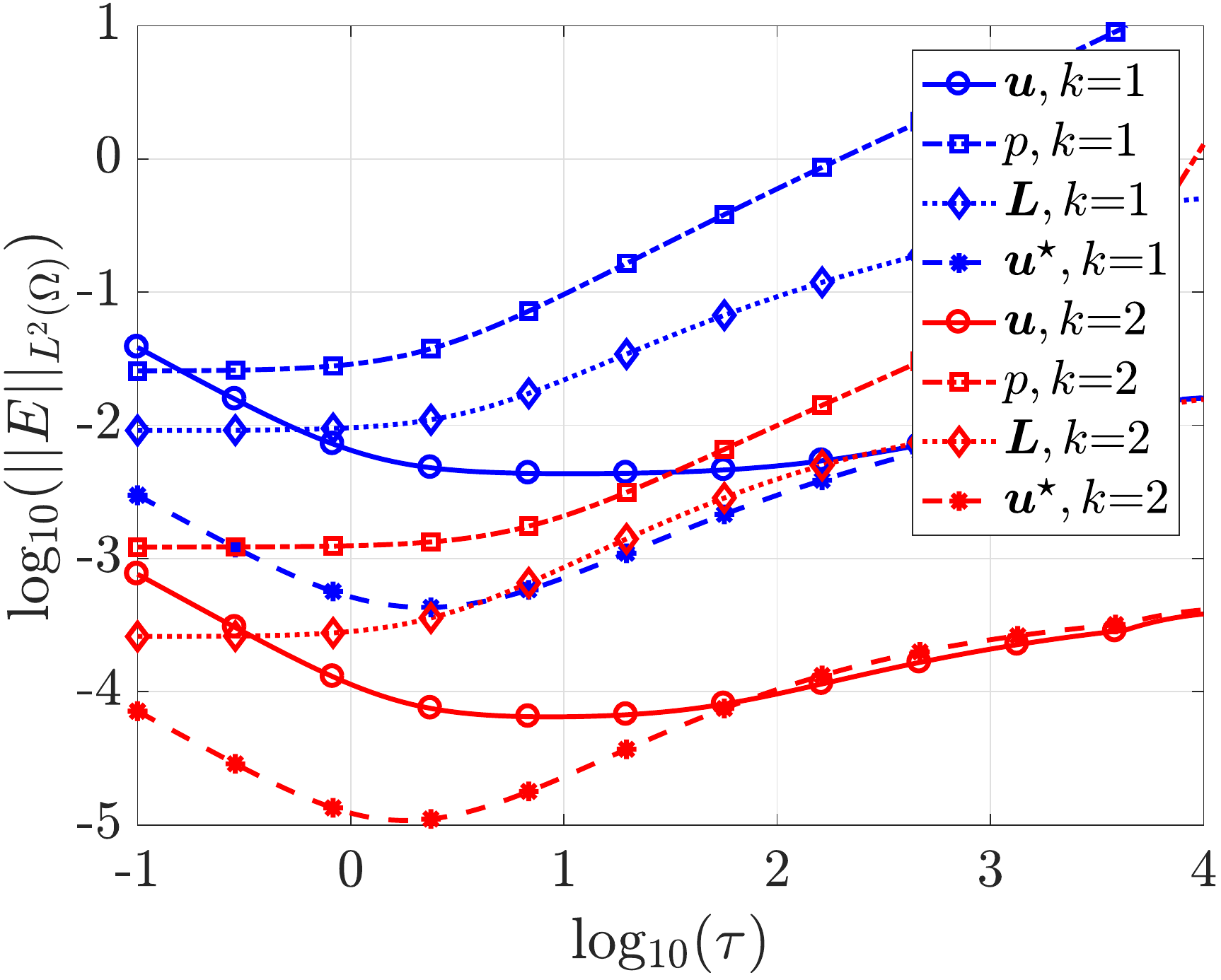}}
	
	\subfigure[Prismatic mesh]{\includegraphics[width=0.4\textwidth]{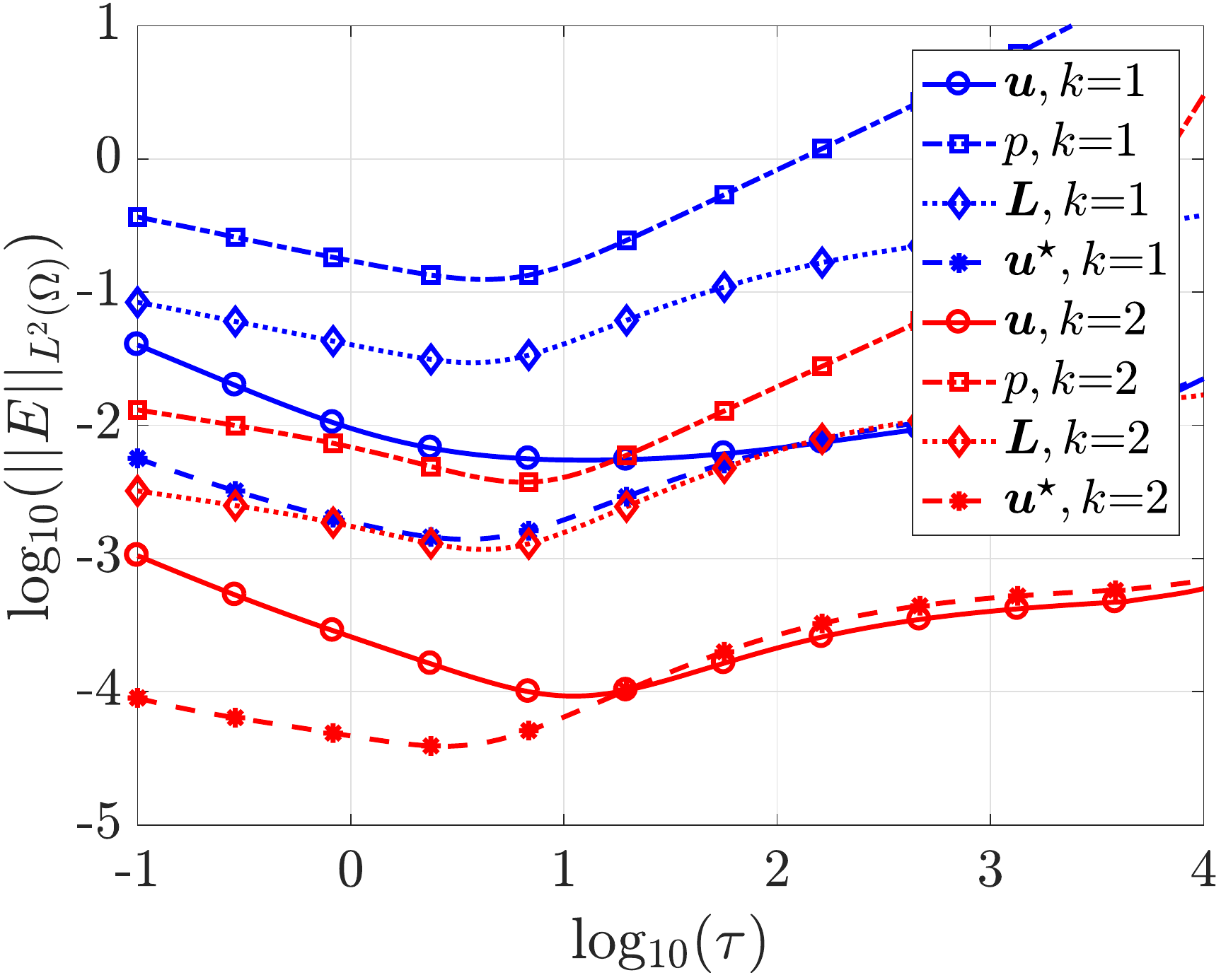}}
	\subfigure[Pyramidal mesh]{\includegraphics[width=0.4\textwidth]{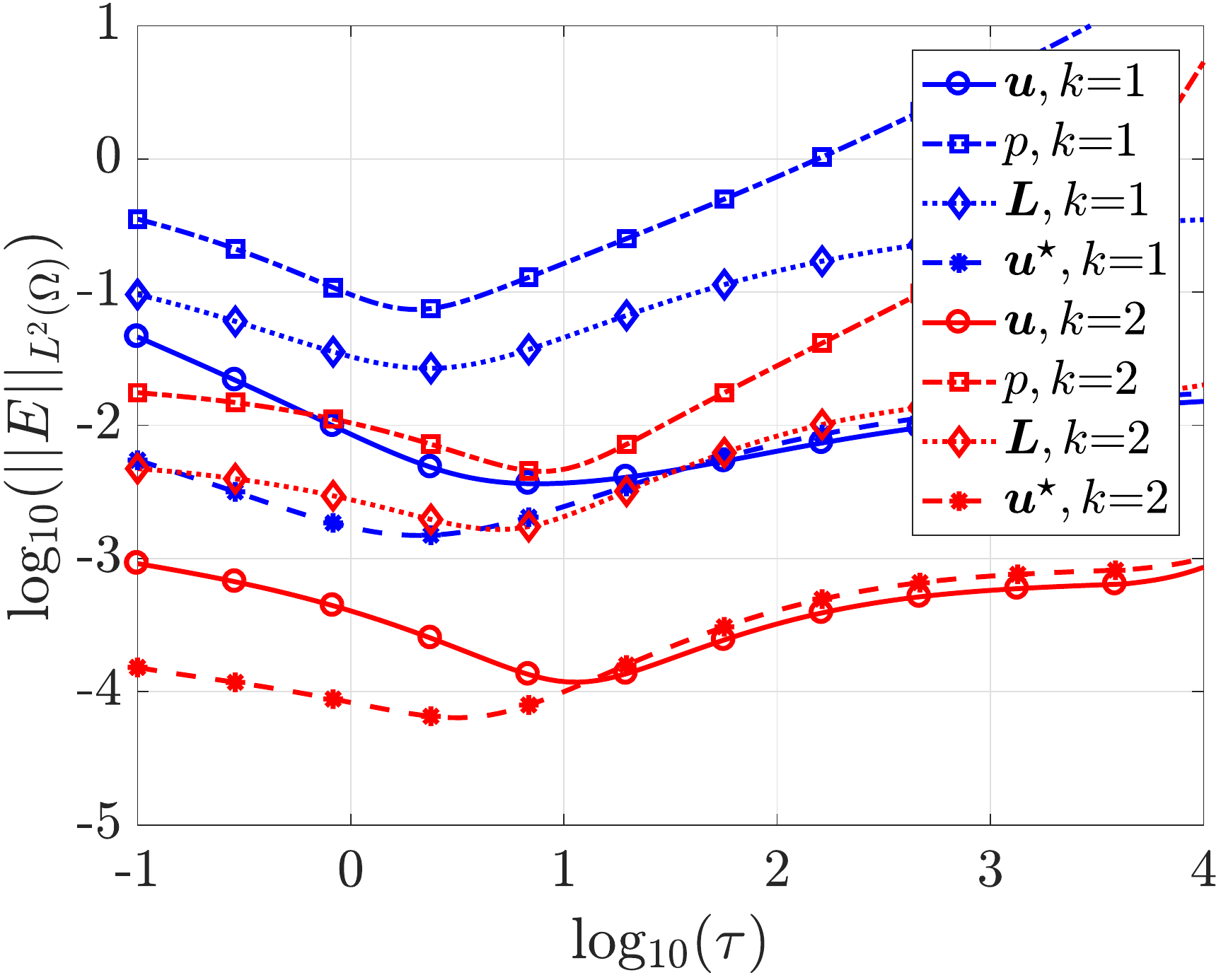}}
	
	\caption{Three dimensional problem: error of the primal, mixed and post-processed variables, $\bu$, $p$, $\bL$ and $\bu^\star$, in the $\eltwo(\Omega)$ norm as a function of the stabilization parameter and for the third level of mesh refinement.}
	\label{fig:tauInfluence3D}
\end{figure}

The evolution of the error of the primal, mixed and post-processed variables, $\bu$, $p$, $\bL$ and $\bu^\star$, in the $\eltwo(\Omega)$ norm as a function of the stabilization parameter $\tau$ is presented in Figure~\ref{fig:tauInfluence3D}.
As highlighted by the theory~\cite{MR2772094} and confirmed by the analysis of the two dimensional case, a value of the stabilization parameter of order one (i.e. $\tau \in [1,10)$) guarantees stability and convergence of the HDG method.
More precisely, a value near $\tau=10$ provides the minimum error for the primal variable but limited or no extra gain in accuracy is obtained through the post-process of the velocity field.
Thus, a value of $\tau=4$ is selected for the following simulations.

The discussed numerical results show that the HDG discretization is robust to the choice of the stabilization parameter.
Moreover, the optimal value of $\tau$ is not dependent upon the degree of approximation or the dimensionality of the problem.
Considering the different types of elements under analysis, the triangular meshes \#2 require a slightly larger value of the stabilization parameter to enter the asymptotic regime and show the optimal convergence and superconvergence properties expected from the theory.

\subsection{Optimal convergence and superconvergence of the primal, mixed and post-processed variables}
\label{sc:convergence}

Consider the optimal values of $\tau$ identified in the previous section.
The optimal convergence properties of the velocity $\bu$, the pressure $p$ and the mixed variable $\bL$ representing the strain rate tensor, are tested for different element types using the $\eltwo(\Omega)$ norm.
Moreover, the superconvergence of the post-processed velocity field $\bu^\star$ is also analyzed.

\subsubsection{Two dimensional example}
\label{sc:convergence2D}

\begin{figure}[!tb]
	\centering
	\subfigure[Quadrilateral meshes: $p, \bL$]{\includegraphics[width=0.4\textwidth]{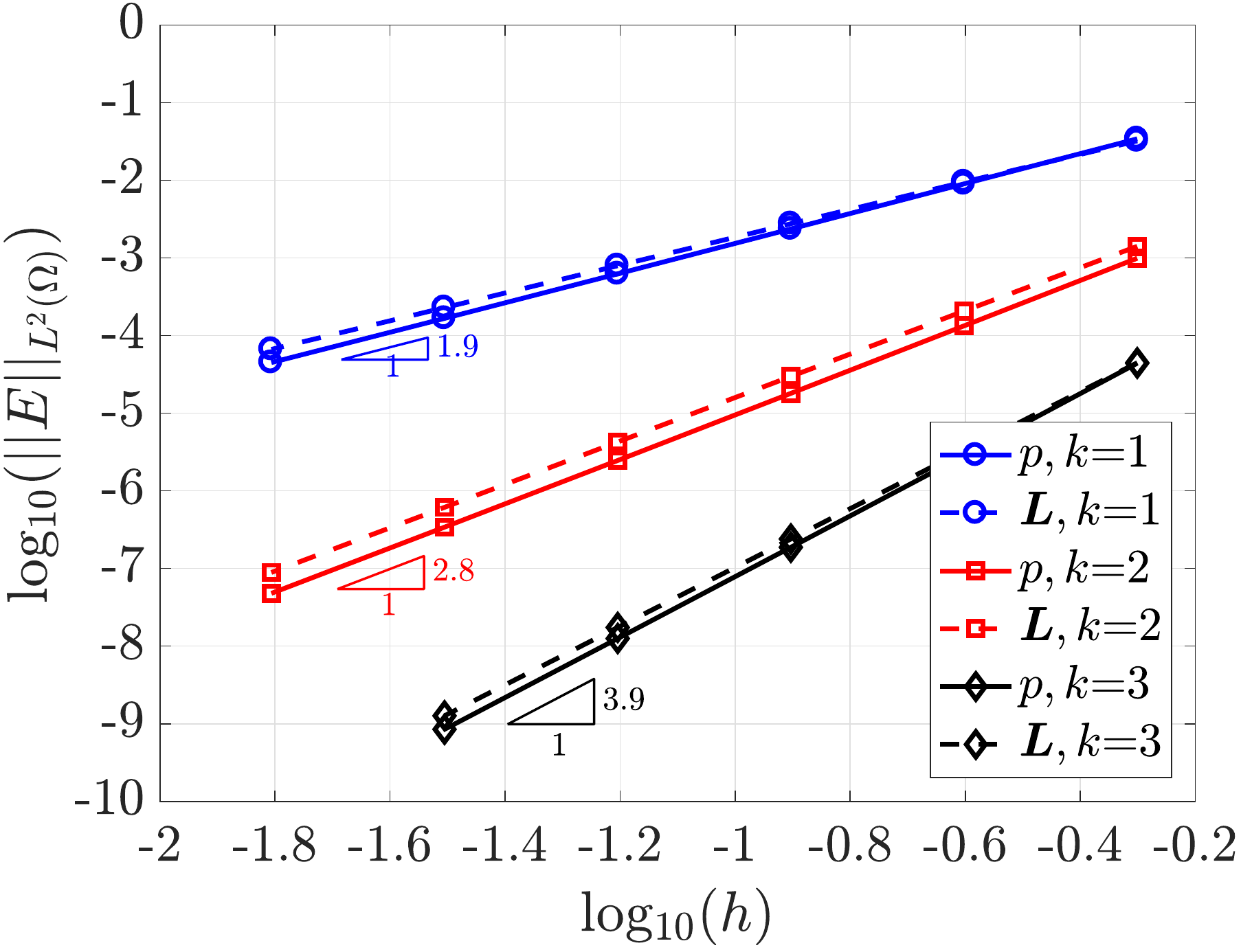}}
	\subfigure[Quadrilateral meshes: $\bu, \bu^\star$]{\includegraphics[width=0.4\textwidth]{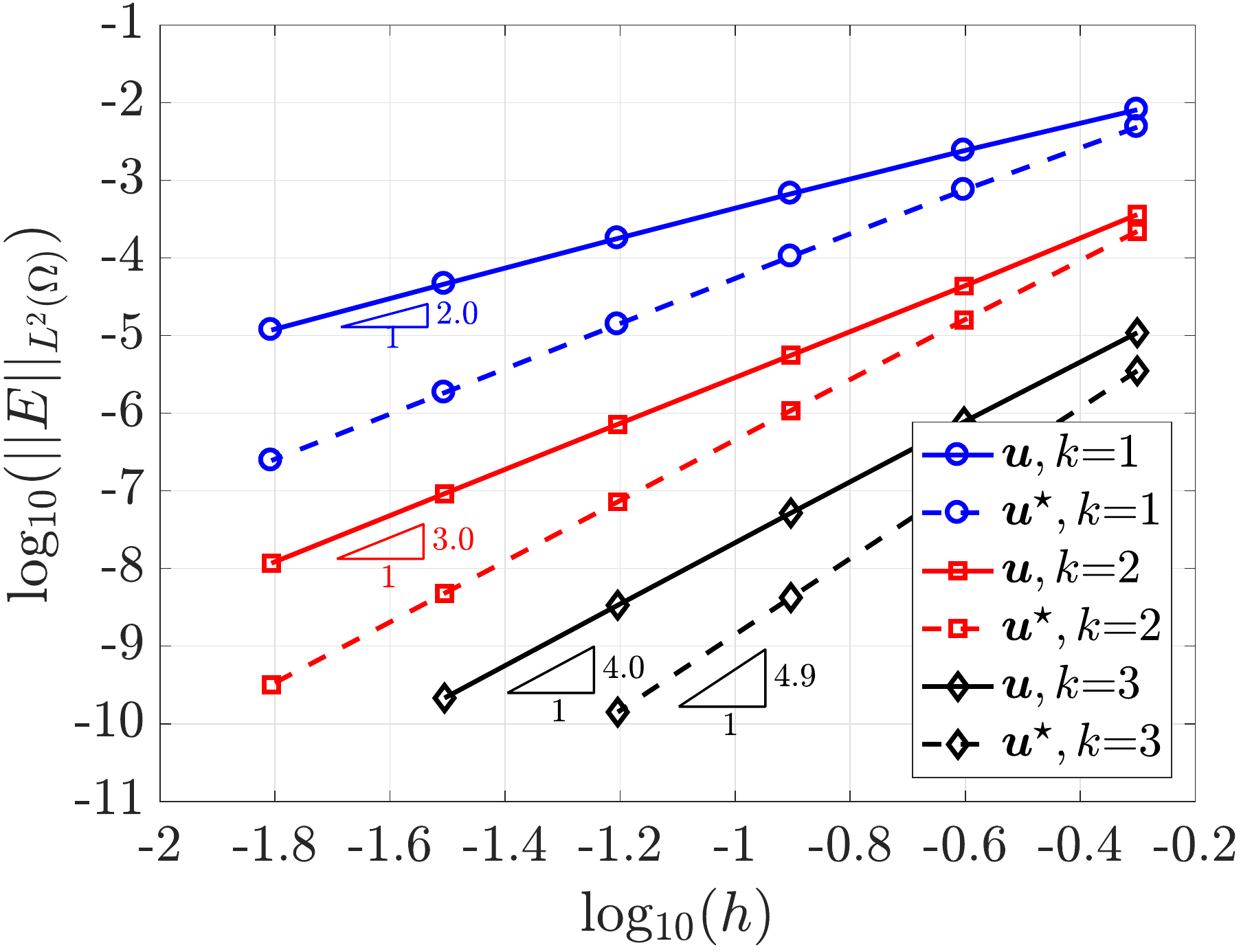}}
	
	\subfigure[Triangular meshes \#1: $p, \bL$]{\includegraphics[width=0.4\textwidth]{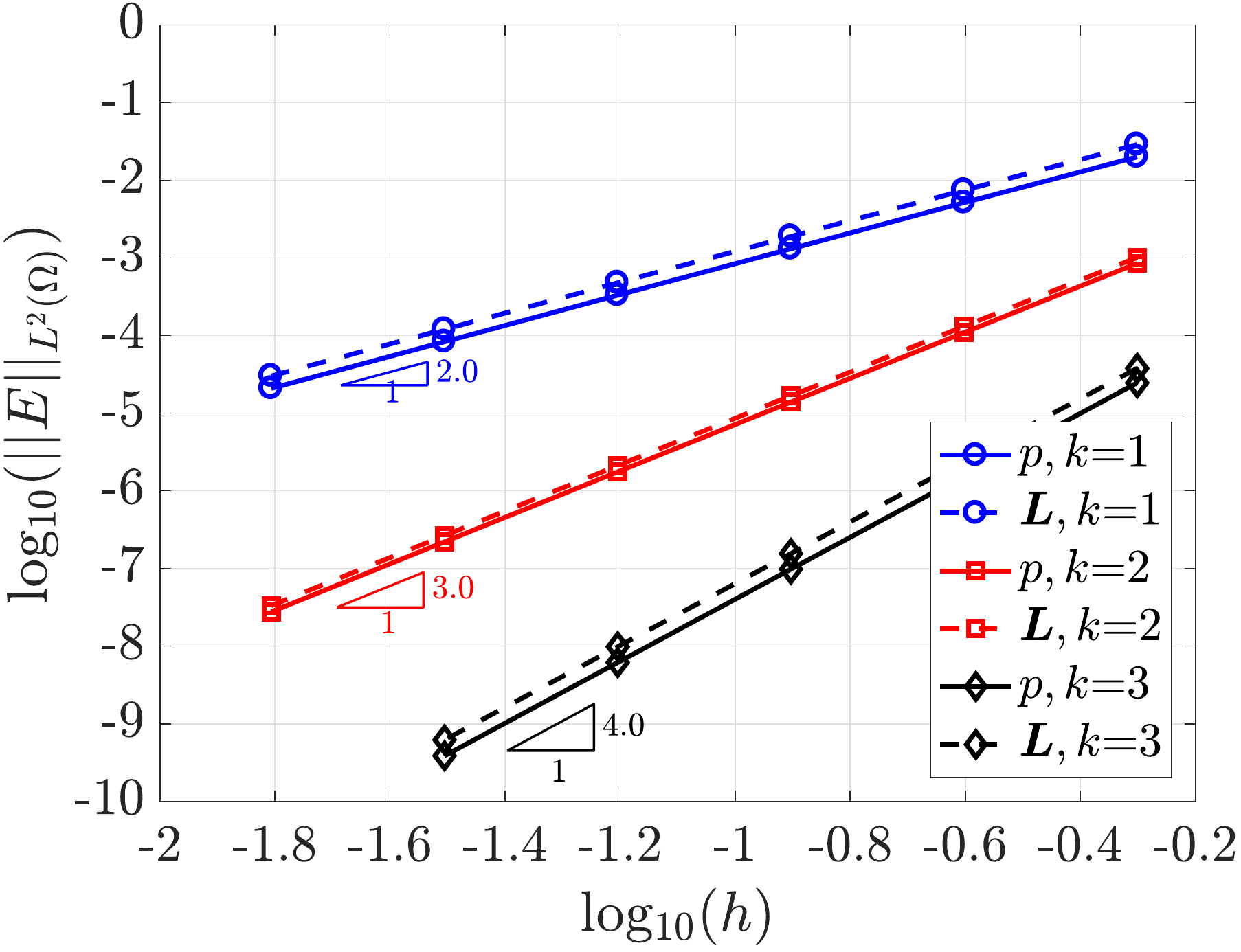}}	
	\subfigure[Triangular meshes \#1: $\bu, \bu^\star$]{\includegraphics[width=0.4\textwidth]{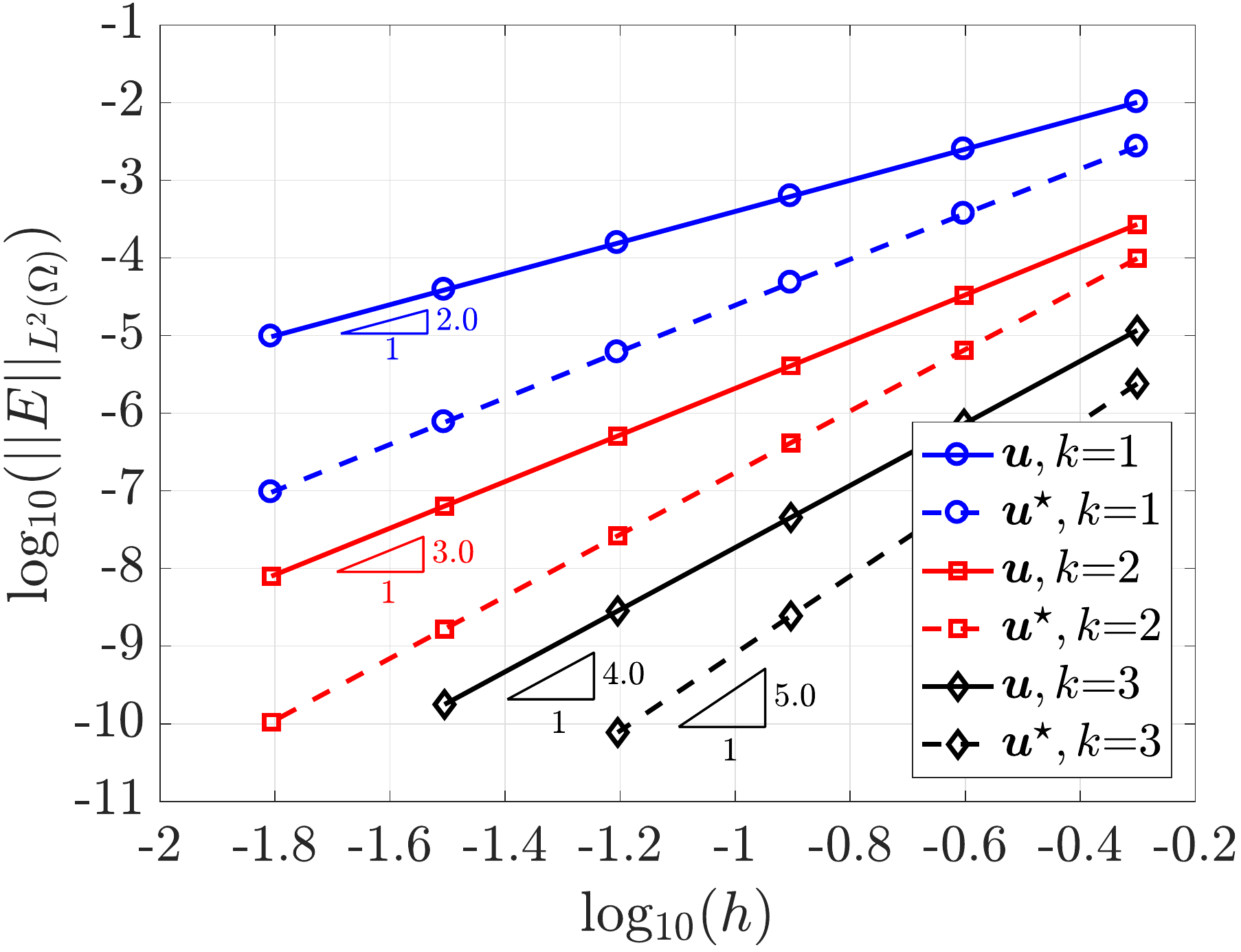}}	
	
	\subfigure[Triangular meshes \#2: $p, \bL$]{\includegraphics[width=0.4\textwidth]{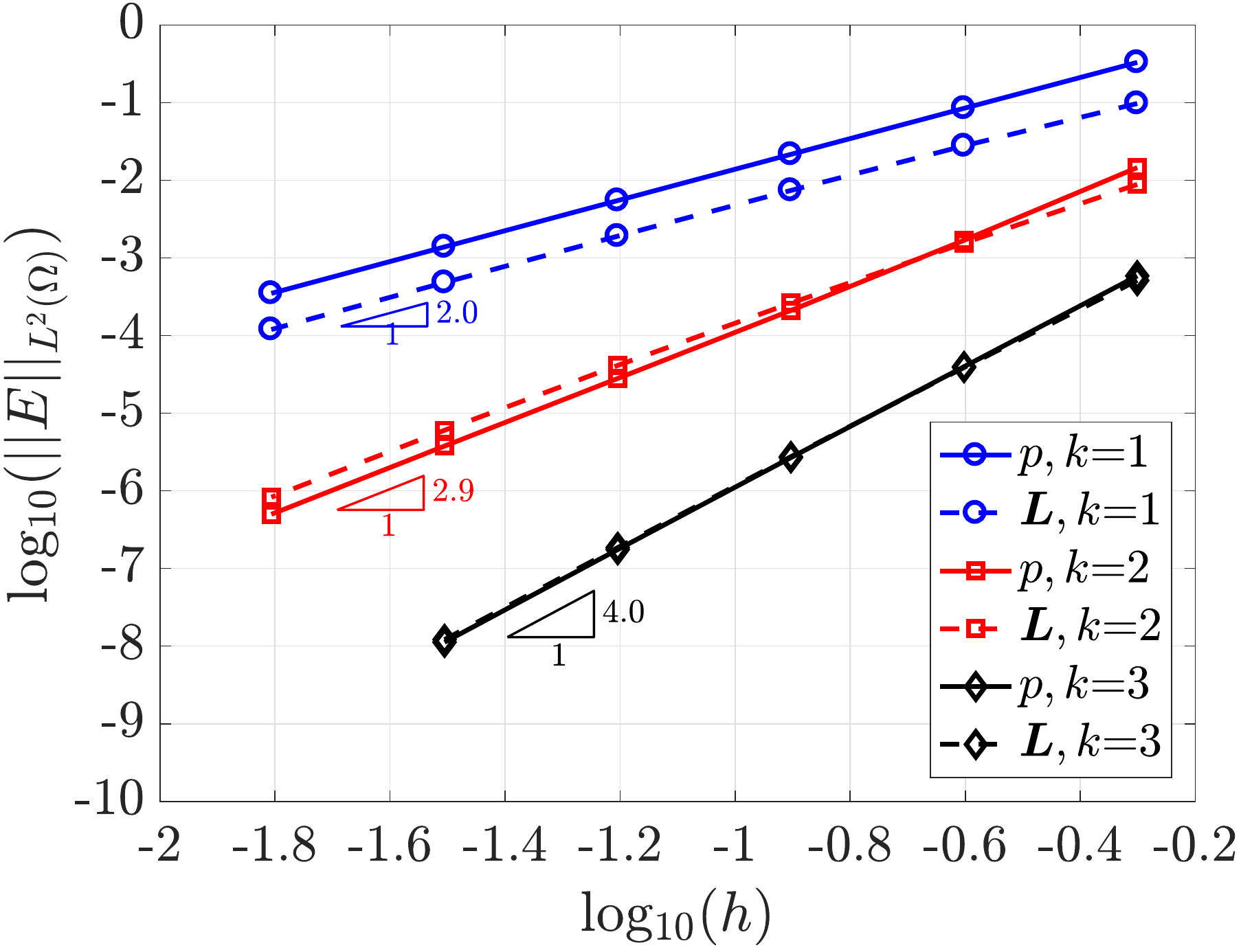}}	
	\subfigure[Triangular meshes \#2: $\bu, \bu^\star$]{\includegraphics[width=0.4\textwidth]{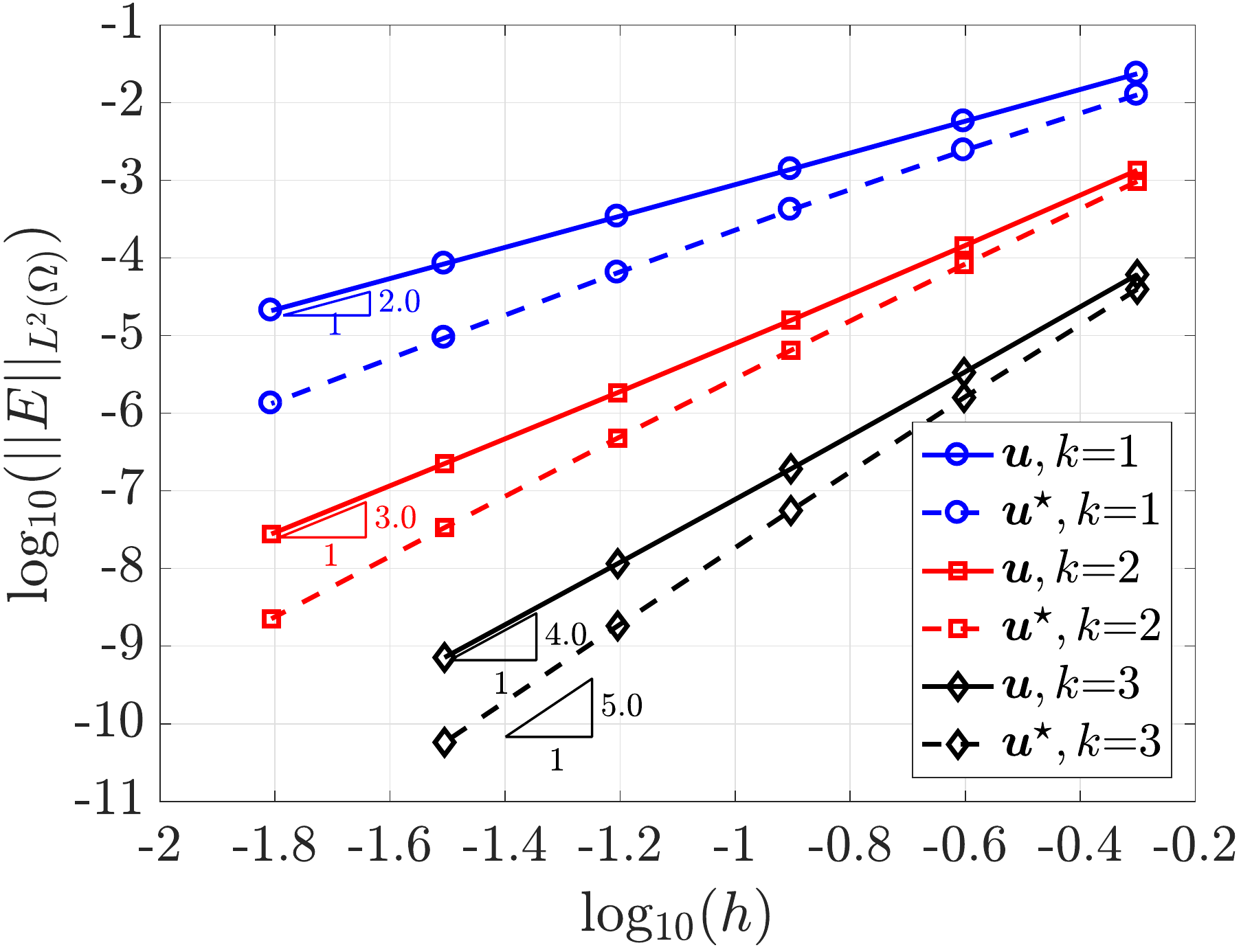}}	
	
	\caption{Two dimensional problem: $h$-convergence of the error of the primal, mixed and post-processed variables, $p$ and $\bL$ (on the left), $\bu$ and $\bu^\star$ (on the right), in the $\eltwo(\Omega)$ norm for quadrilateral and triangular meshes with different degrees of approximation.}
	\label{fig:hConv2D}
\end{figure}

In Figure~\ref{fig:hConv2D}, the first column presents the convergence of the error of the primal and mixed variables $p$ and $\bL$, measured in the $\eltwo(\Omega)$ norm, as a function of the characteristic element size $h$ for both quadrilateral and triangular elements and for a degree of approximation ranging from $k=1$ up to $k=3$. 
In a similar fashion, the second column provides the corresponding convergence history for the primal and the post-processed velocities $\bu$ and $\bu^\star$.

It can be observed that almost the optimal or the optimal rate of convergence $h^{k+1}$ is obtained for $\bu$, $p$ and $\bL$, for all the element types and degrees of approximation considered. 
As previously mentioned, the triangular mesh \#1 has considerably more degrees of freedom than the other meshes: in particular, for the same characteristic element size, the triangular meshes \#1 have approximately 2.5 times (respectively, 5 times) more internal faces than the triangular mesh \#2 (respectively, the quadrilateral mesh). 
Thus, despite the results in Figure~\ref{fig:hConv2D} indicate that the triangular mesh \#1 provides more accuracy than the other meshes, a comparison in terms of the global number of degrees of freedom confirms that similar results are obtained using meshes of different element types.
Concerning the post-processed variable, the rate of convergence $h^{k+2}$ is achieved and the superconvergence property is verified.
This confirms that the average of the hybrid variable $\bhu$ on the boundary leads to a superconvergent approximation, as observed in \cite{preprintVoigtElasticity} for the linear elastic problem.
Beside the improved convergence rate, the discussed post-process procedure is responsible for a gain in accuracy of $\bu^\star$ with respect to the original approximation $\bu$ of the velocity field.
Hence, the information encapsulated in the primal and post-processed variables may be exploited to construct an error indicator and devise an automatic degree adaptivity strategy as discussed in \cite{giorgiani2014hybridizable,RS-SH:18}.

\subsubsection{Three dimensional example}
\label{sc:convergence3D}

\begin{figure}[!tb]
	\centering
	\subfigure[Hexahedral meshes: $p, \bL$]{\includegraphics[width=0.4\textwidth]{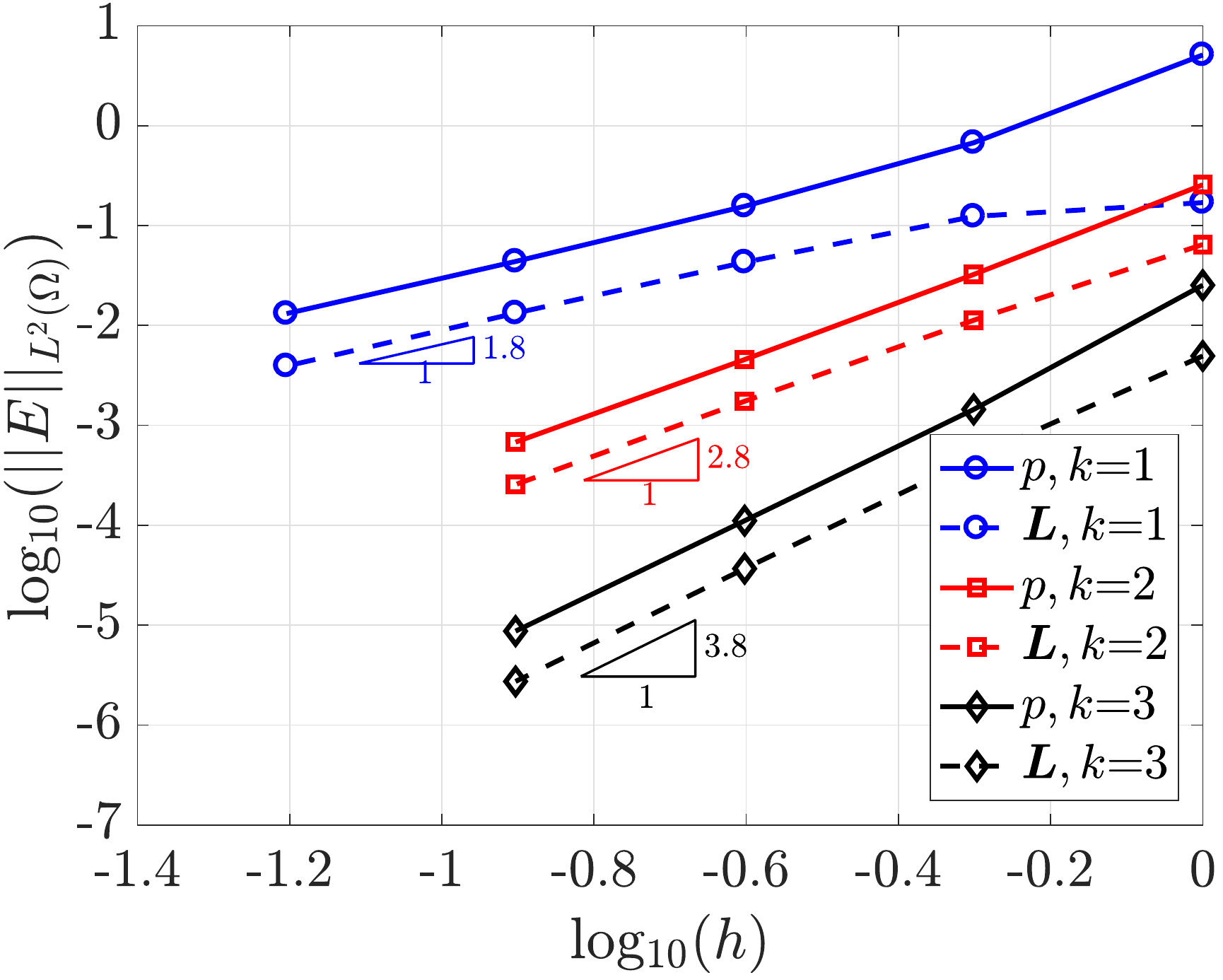}}
	\subfigure[Tetrahedral meshes: $p, \bL$]{\includegraphics[width=0.4\textwidth]{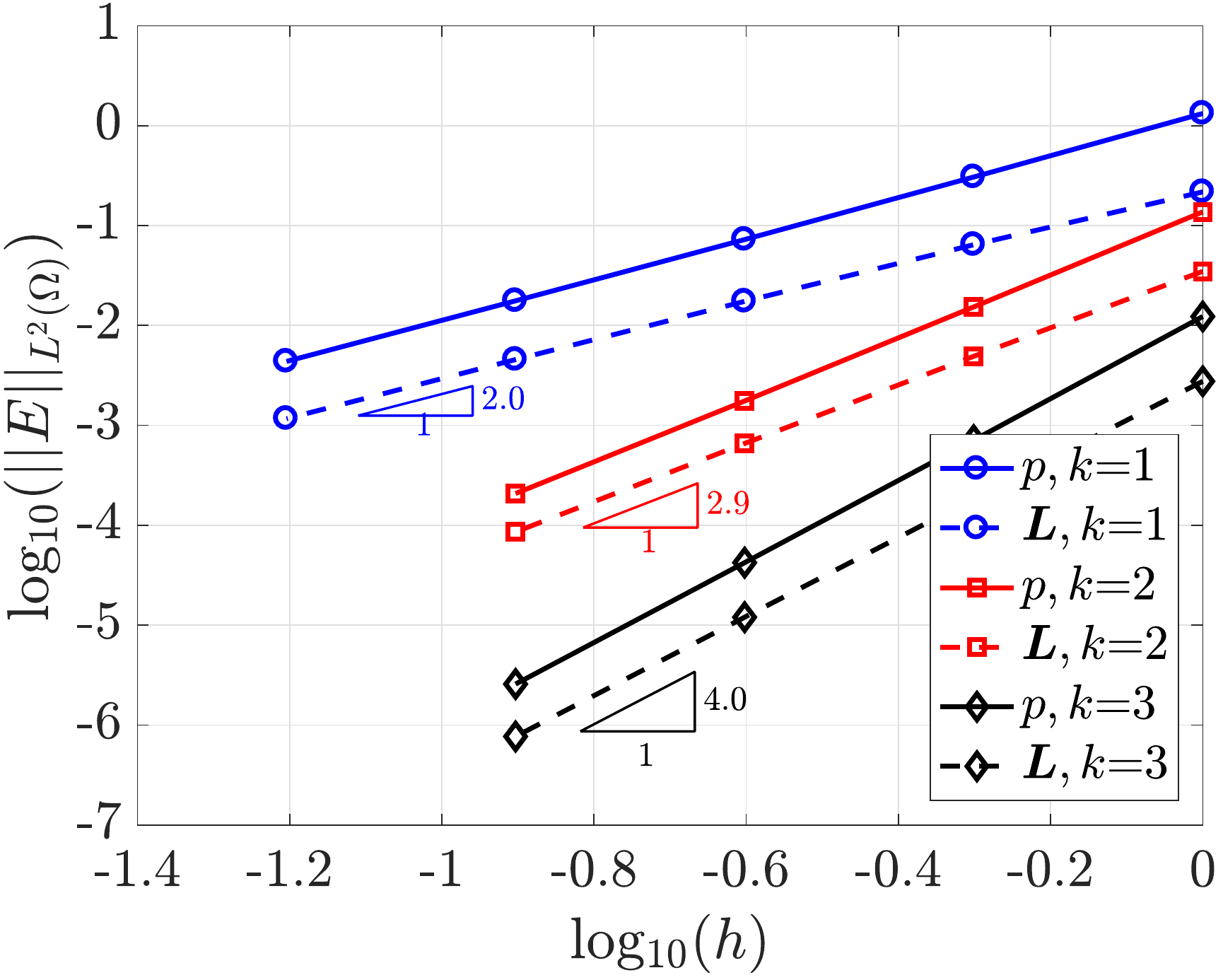}}
	
	\subfigure[Prismatic meshes: $p, \bL$]{\includegraphics[width=0.4\textwidth]{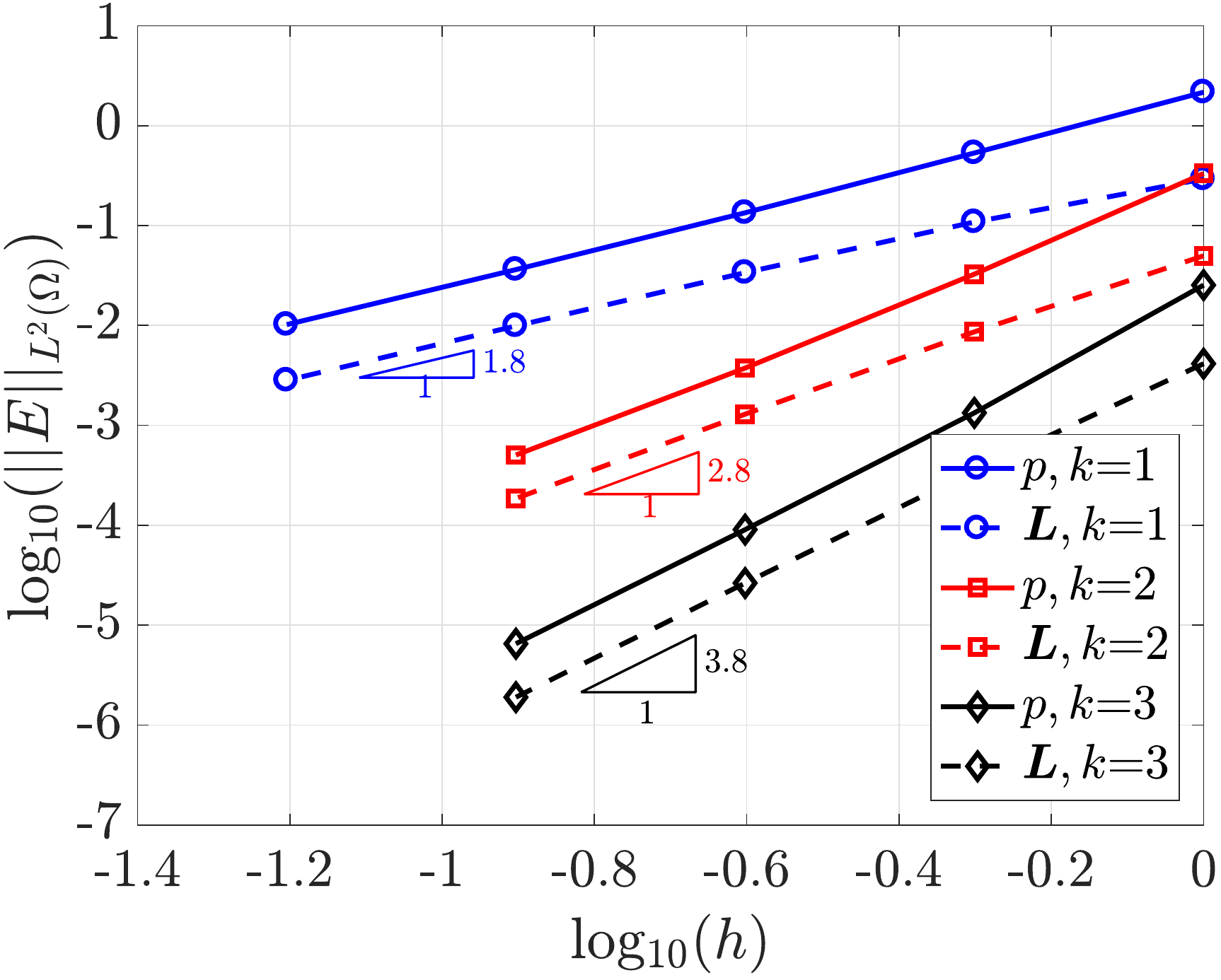}}
	\subfigure[Pyramidal meshes: $p, \bL$]{\includegraphics[width=0.4\textwidth]{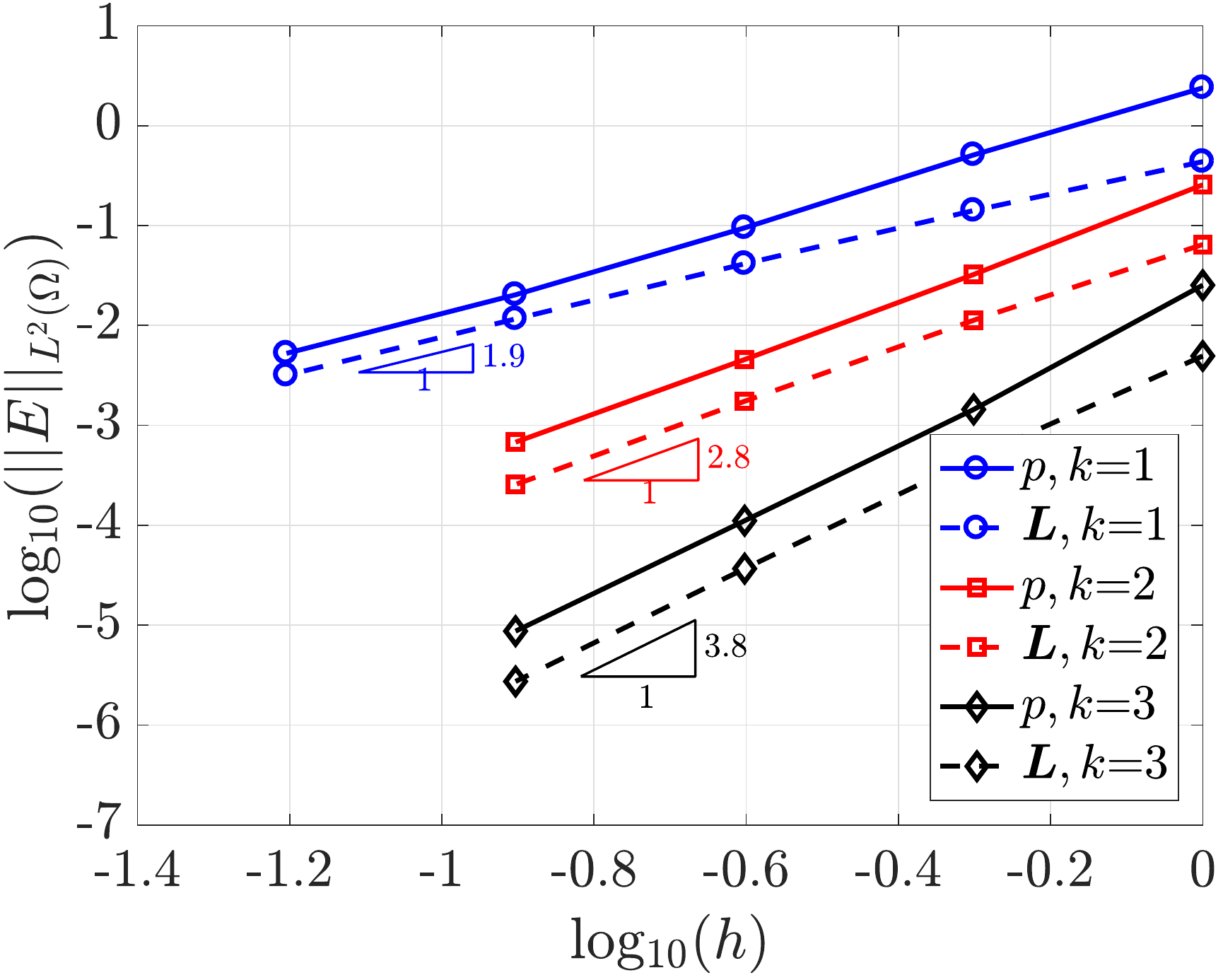}}
	
	\caption{Three dimensional problem: $h$-convergence of the error of the primal and mixed variables, $p$ and $\bL$, in the $\eltwo(\Omega)$ norm for hexahedral, tetrahedral, prismatic and pyramidal  meshes with different degrees of approximation.}
	\label{fig:hConv3D_pL}
\end{figure}

\begin{figure}[!tb]
	\centering
	\subfigure[Hexahedral meshes: $\bu, \bu^\star$]{\includegraphics[width=0.4\textwidth]{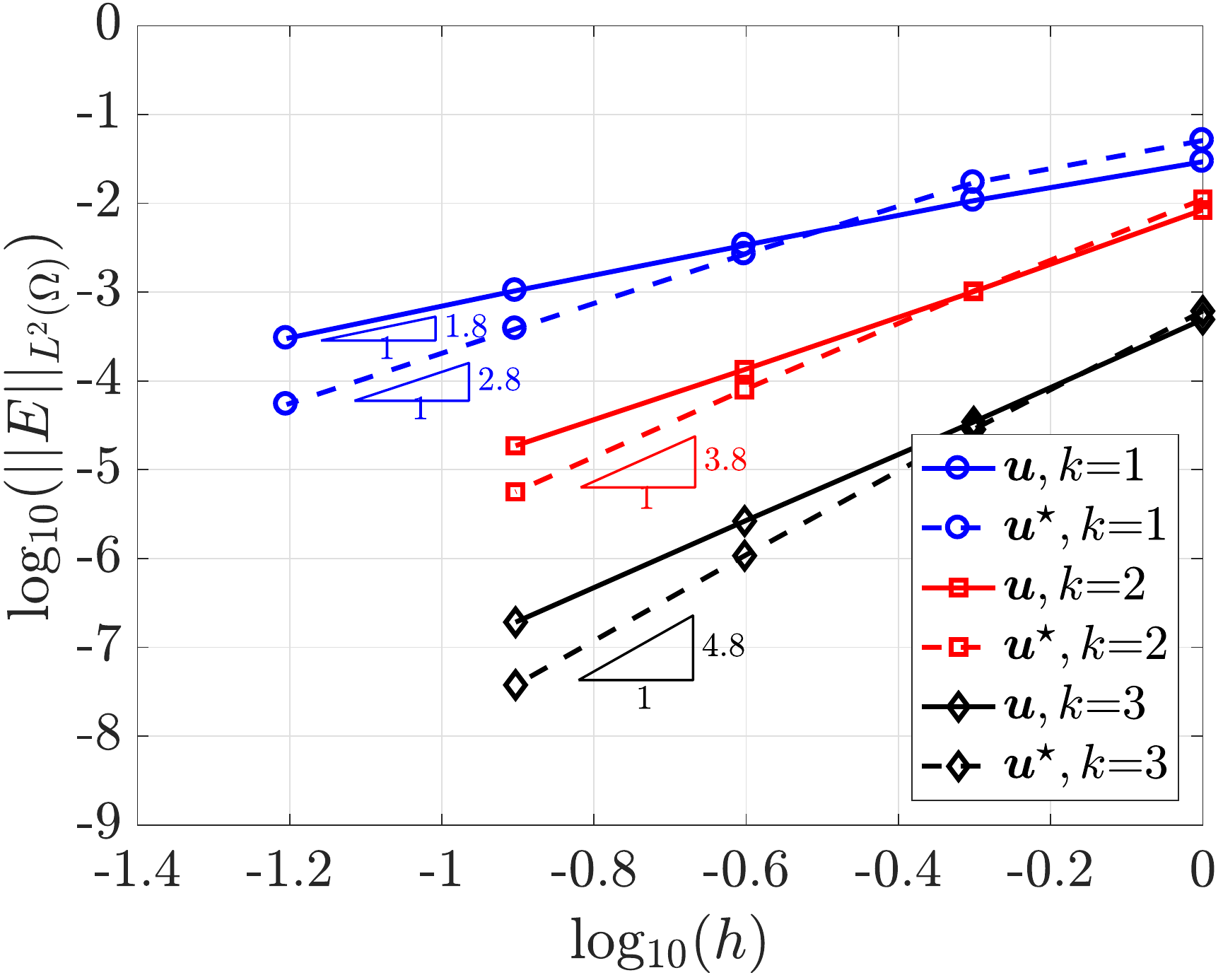}}
	\subfigure[Tetrahedral meshes: $\bu, \bu^\star$]{\includegraphics[width=0.4\textwidth]{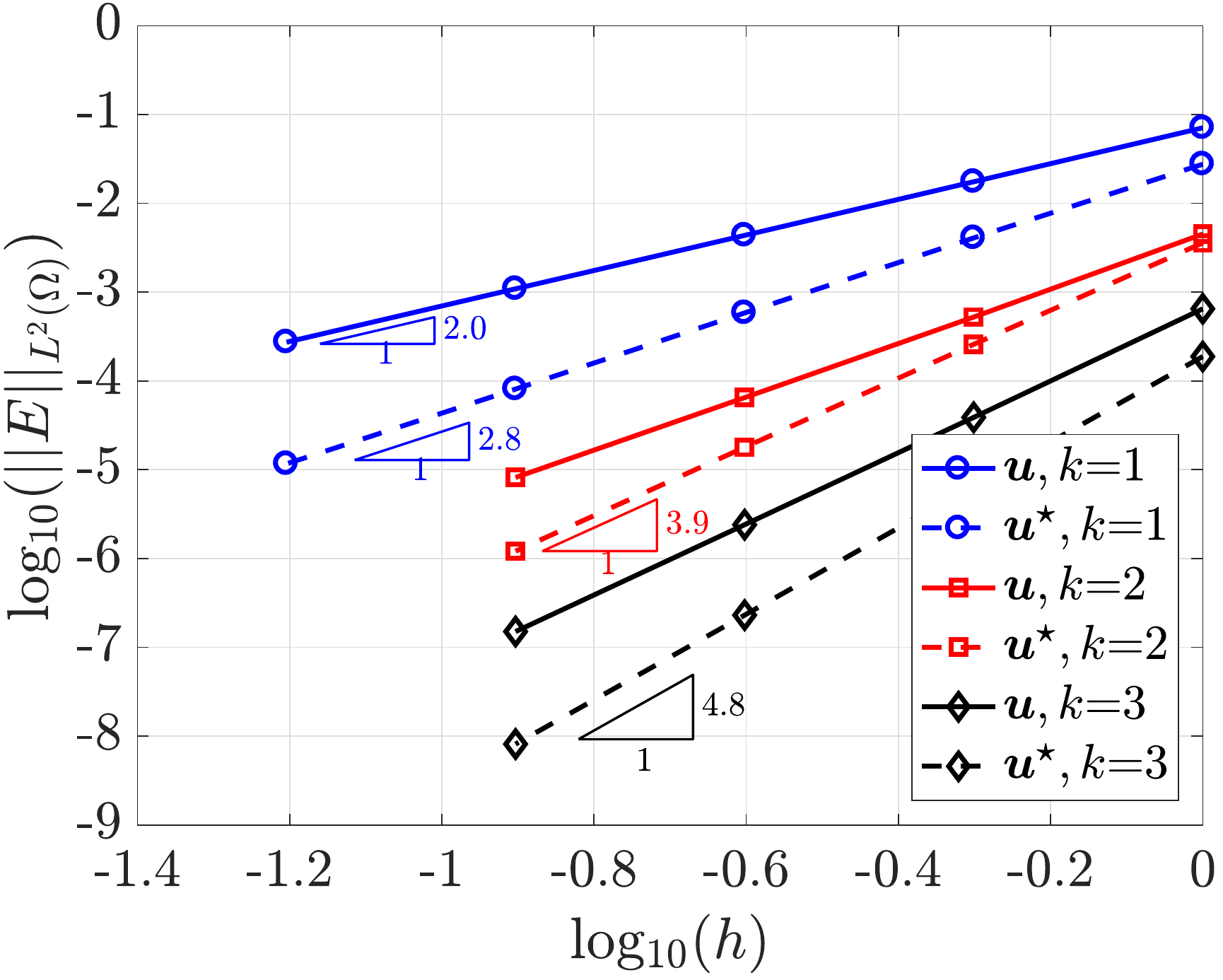}}
	
	\subfigure[Prismatic meshes: $\bu, \bu^\star$]{\includegraphics[width=0.4\textwidth]{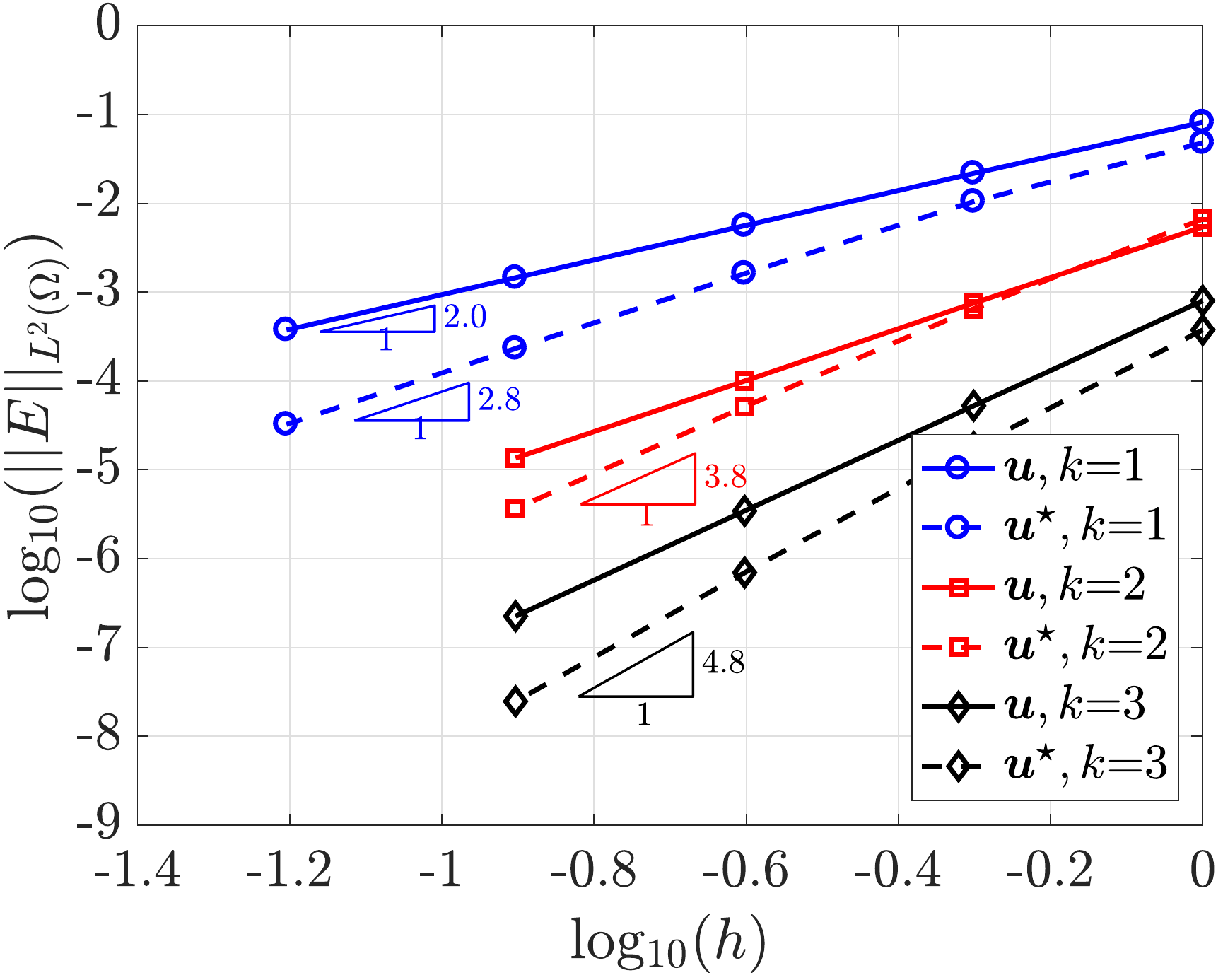}}
	\subfigure[Pyramidal meshes: $\bu, \bu^\star$]{\includegraphics[width=0.4\textwidth]{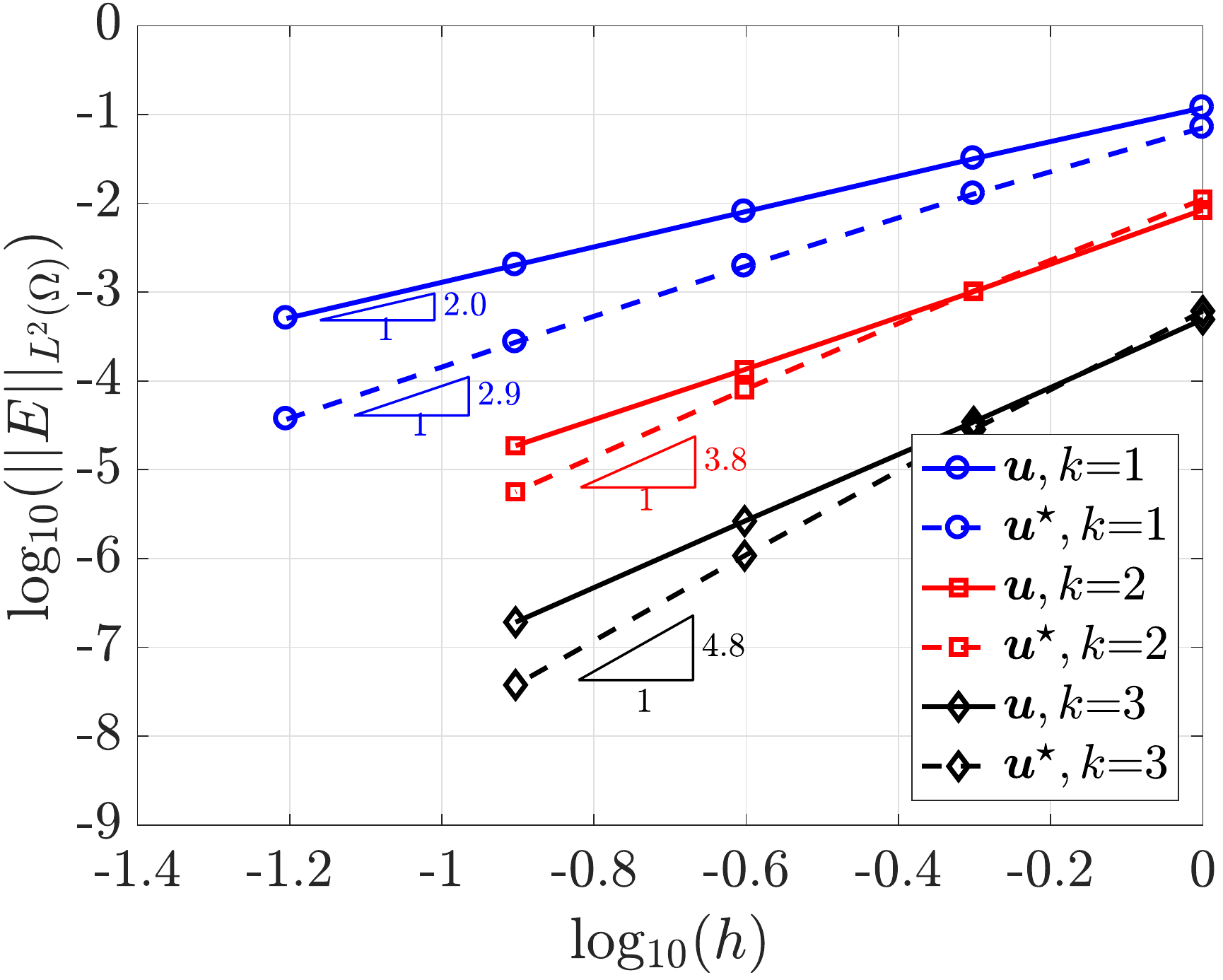}}
	
\caption{Three dimensional problem: $h$-convergence of the error of the primal and post-processed variables, $\bu$ and $\bu^\star$, in the $\eltwo(\Omega)$ norm for hexahedral, tetrahedral, prismatic and pyramidal meshes with different degrees of approximation.}
	\label{fig:hConv3D_uU}
\end{figure}

Similarly to the previous example, the convergence of the error of $p$ and $\bL$ (Fig.~\ref{fig:hConv3D_pL}) and $\bu$ and $\bu^\star$ (Fig.~\ref{fig:hConv3D_uU}), measured in the $\eltwo(\Omega)$ norm, as a function of the characteristic element size $h$ is presented for hexahedral, tetrahedral, prismatic and pyramidal elements and for a degree of approximation ranging from $k=1$ up to $k=3$. 

As for the two dimensional case, almost the optimal or the optimal rate of convergence $h^{k+1}$ is obtained for $\bu$, $p$ and $\bL$ in 3D, for all the element types and degrees of approximation considered (cf. Fig.~\ref{fig:hConv3D_pL}-\ref{fig:hConv3D_uU}). 
In Figure~\ref{fig:hConv3D_uU}, the post-processed variable is shown to superconverge with a rate of convergence $h^{k+2}$.
Beside the improved convergence rate, the discussed post-process procedure is responsible of a gain in accuracy of $\bu^\star$ with respect to the original approximation $\bu$ of the velocity field.

The presented numerical experiments in two and three dimensions confirm that exploiting Voigt notation the HDG approximation of the Stokes equation achieves optimal convergence rate $h^{k+1}$ for both the primal variables $\bu$ and $p$ and the mixed one $\bL$. 
In particular, contrary to what observed in \cite{Nguyen-CNP:10}, the convergence of the mixed variable does not deteriorate when considering the Cauchy formulation of the Stokes flow.
As discussed in \cite{preprintVoigtElasticity} for the linear elastic problem, the post-process technique exploiting the $\curl$ of $\bu$ allows to construct an approximation of the primal vector field superconverging with order $k+2$.
Moreover, the post-process strategy provides an extra gain in accuracy with respect to the original approximation of the velocity field.
As highlighted in Figure~\ref{fig:hConv3D_uU}, a solution that is almost one order of magnitude more precise than the HDG solution is obtained, even for linear approximations.

\subsection{Numerical evaluation of quantities of interest: drag force on a sphere}
\label{sc:sphere}

The last example considers the classical test case of the viscous flow around a sphere. 
The objective of this test is to show the capability of the described HDG method to provide an approximation of the pressure and the viscous forces sufficiently accurate to evaluate a quantity of interest with the precision required by industrial standards.
Consider the domain $\Omega = \left( [-7,15] \times [-5,5] \times [-5,5] \right)\setminus \mathcal{B}_{1,\bm{0}}$, $\mathcal{B}_{1,\bm{0}}$ being a ball of unit radius centered at the origin.
\begin{figure}[!tb]
	\centering
	\subfigure[Magnitude of the velocity with streamlines]{\includegraphics[width=0.4\textwidth]{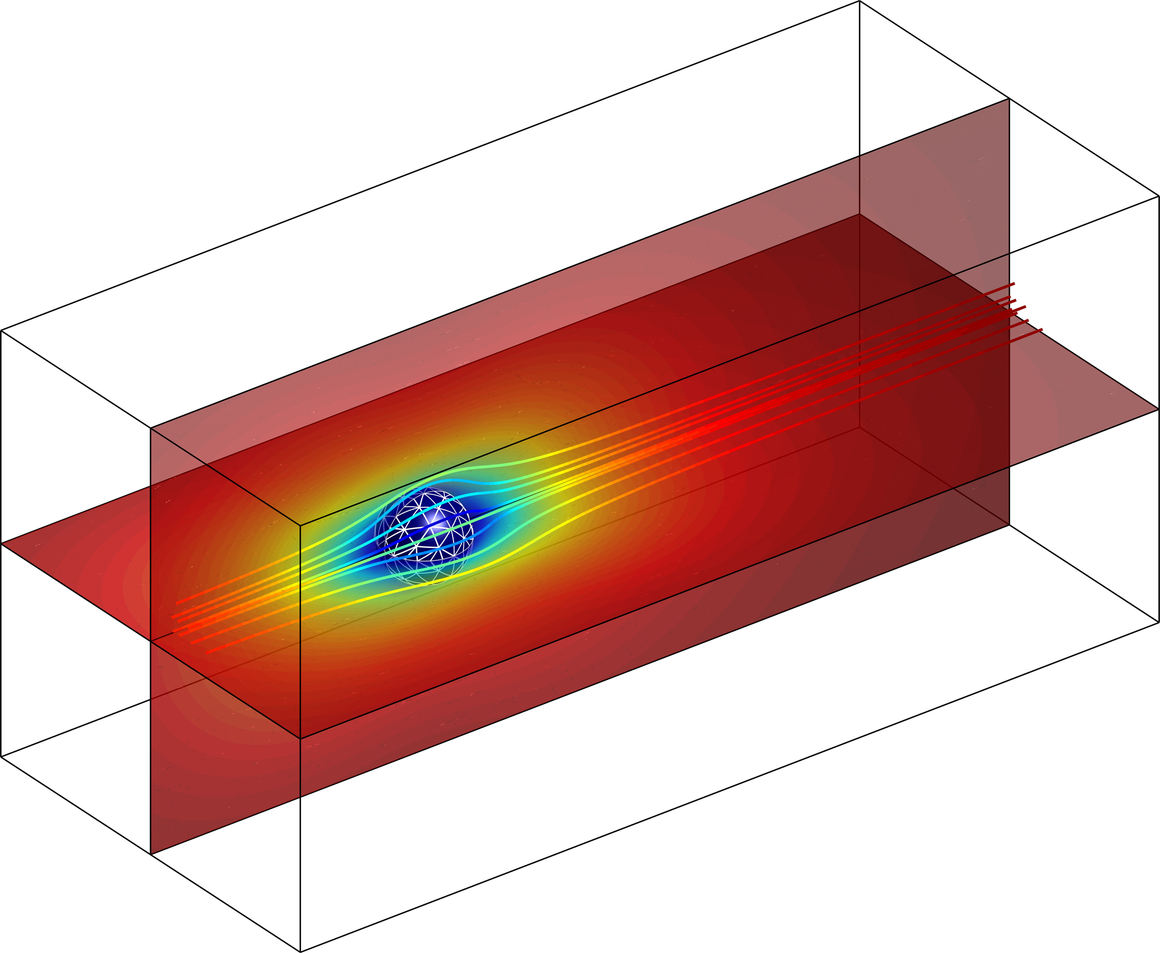}}
	\subfigure[Pressure field]{\includegraphics[width=0.4\textwidth]{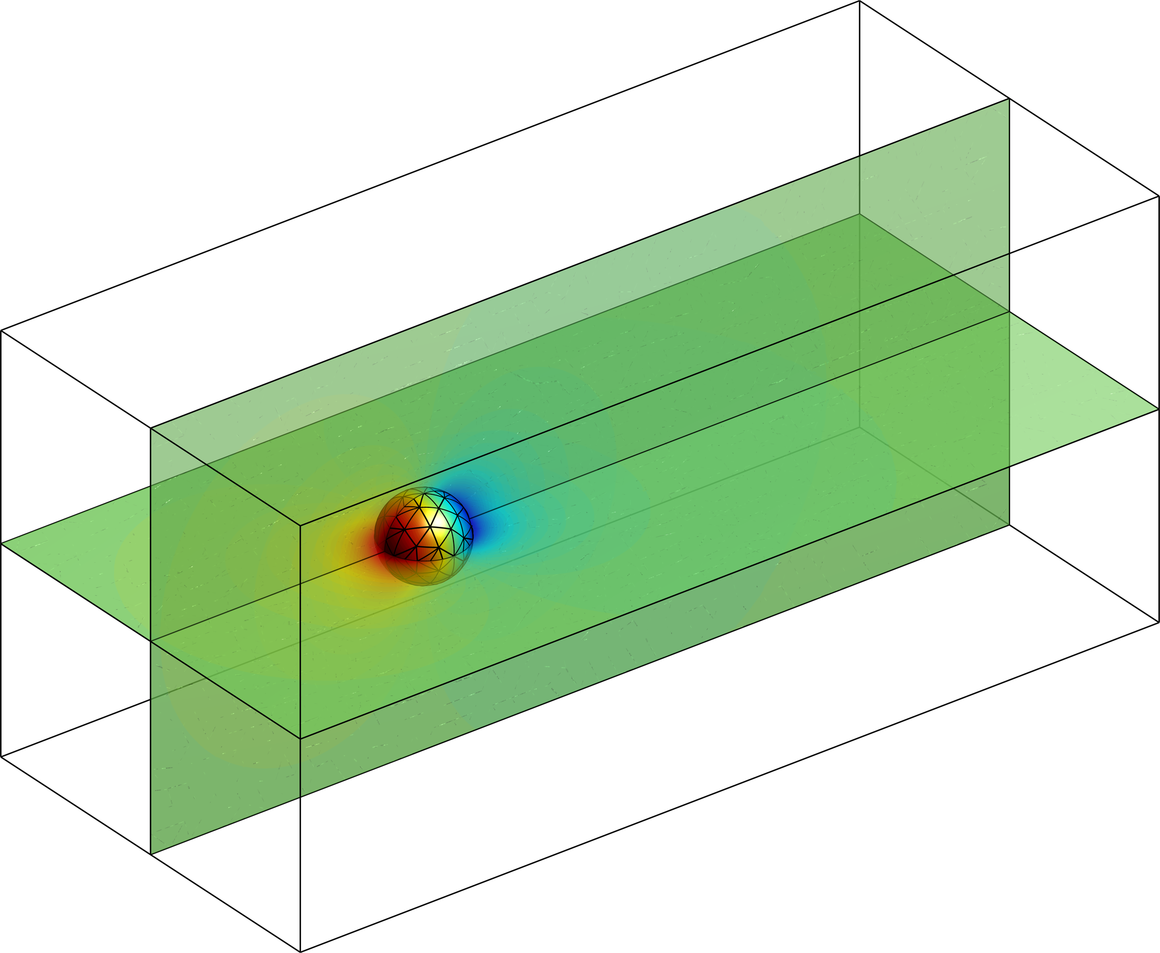}}
	
	\subfigure[Drag force]{\includegraphics[width=0.4\textwidth]{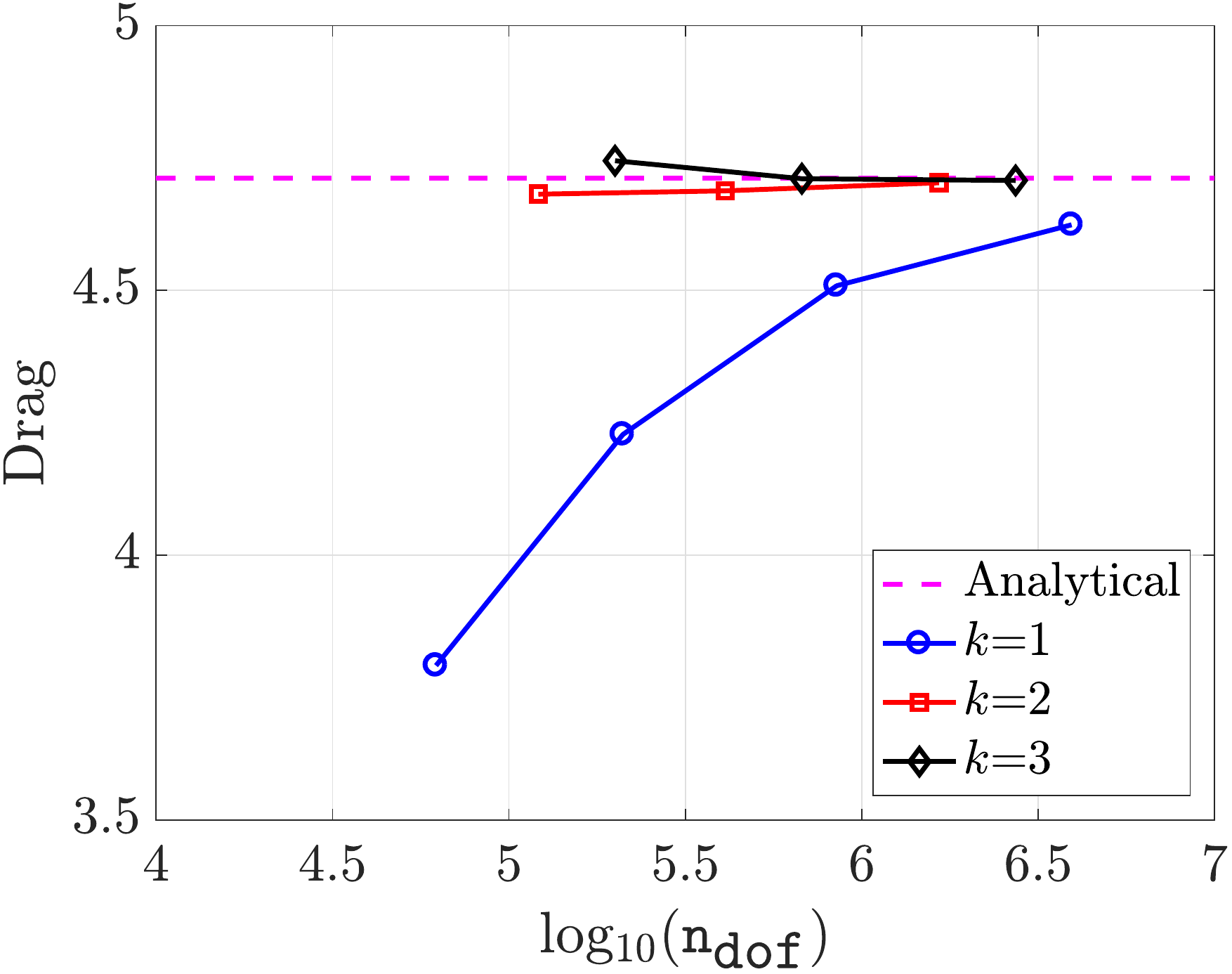}}
	\caption{Flow past a sphere: HDG approximation of (a) the velocity field with streamlines of the flow and (b) the pressure field using the third level of refinement of a tetrahedral mesh and $k=2$. 
	(c) Convergence of the drag as a function of the number of degrees of freedom.}
	\label{fig:3DsphereSol}
\end{figure}
To reduce the computational effort, the symmetry of $\Omega$ is exploited and solely one fourth of the domain is taken into account to perform the numerical experiments.
Different tetrahedral meshes of the domain are considered, ranging from 3,107 to 204,099 elements.
High-order computations employ isoparametric curved meshes. The extension to high-order is performed using the solid mechanics analogy described in \cite{poya2016unified,xie2013generation}.
Figure~\ref{fig:3DsphereSol} (a)-(b) shows the magnitude of the velocity with streamlines of the flow and the pressure field computed on the third level of refinement of the mesh, featuring 43,682 tetrahedrons, and using a quadratic degree of approximation.

The results in Figure~\ref{fig:3DsphereSol} (c) show the convergence of the drag force as the number of degrees of freedom is increased, i.e. for different levels of mesh refinement and for a degree of approximation ranging from $k=1$ up to $k=3$. 
The numerically computed drag is compared with the analytical value from the literature~\cite{batchelor2000introduction}.
In Table~\ref{tab:dragSphere}, a quantitative analysis of the relative error in the computation of the drag force is reported for all the mesh refinements and degrees of approximation considered. 
Using linear elements, almost 4 millions degrees of freedom are required by the method to compute the drag coefficient with a relative error of 2\%.
The same level of accuracy is achieved by quadratic and cubic elements using the coarsest mesh under analysis and less than 200,000 degrees of freedom.
More precisely, moving to high-order approximations, errors lower than 0.5\% are obtained using few hundreds thousands degrees of freedom.
The observed additional accuracy results from the concurrent use of high-order polynomial functions for the discretization of the unknown variables and high-order approximations of the geometry via meshes featuring curved elements.
Thus, the superiority of high-order methods with respect to low-order ones discussed in the literature (cf. e.g. \cite{SEVILLA201315}) is confirmed.
\begin{table}[hbt]
\centering
\begin{tabular}[hbt]{| c | c | c | c | c |}
	\hline
	$k$ & Mesh & Elements & $\ndof$ & Drag error \\
	\hline
	\multirow{4}{*}{1} & 1 & 3,107 & 62,147 & $1.95 \cdot 10^{-1}$ \\
	\cline{2-5}
	& 2 & 10,680 & 210,453 & $1.03 \cdot 10^{-1}$ \\
	\cline{2-5}
	& 3 & 43,682 & 849,452 &  $4.32 \cdot 10^{-2}$ \\
	\cline{2-5}
	& 4 & 204,099 & 3,934,212 & $1.88 \cdot 10^{-2}$ \\
	\hline
	\multirow{3}{*}{2} & 1 & 3,107 & 121,187 & $6.52  \cdot 10^{-3}$ \\
	\cline{2-5}
	& 2 & 10,680 & 410,226 &  $5.18 \cdot 10^{-3}$ \\
	\cline{2-5}
	& 3 & 43,682 & 1,655,222 &  $1.96 \cdot 10^{-3}$ \\
	\hline
	\multirow{3}{*}{3} & 1 & 3,107 & 199,907 &  $6.88 \cdot 10^{-3}$ \\
	\cline{2-5}
	& 2 & 10,680 & 676,590 & $4.25 \cdot 10^{-4}$ \\
	\cline{2-5}
	& 3 & 43,682 & 2,729,582 & $1.02 \cdot 10^{-3}$ \\
	\hline
\end{tabular}
\caption{Flow past a sphere: relative error in the computation of the drag force for different levels of mesh refinement and with different degrees of approximation.}
\label{tab:dragSphere}
\end{table}

\section{Conclusion}
\label{sc:Conclusion}

This paper describes a hybridizable discontinuous Galerkin method using Voigt notation, first introduced in \cite{preprintVoigtElasticity}, for the Cauchy formulation of the Stokes equation. 
Owing to Voigt notation, the symmetry of the stress tensor is strongly enforced by storing in a vector format only half of the off-diagonal terms.
Moreover, physically meaningful tractions may be naturally imposed on the Neumann boundary.
Contrary to the existing superconvergent HDG formulations involving the symmetric part of the gradient, the proposed method does not enrich the discrete spaces of approximation and it reduces the number of degrees of freedom of the mixed variable.
Hence, the resulting local problems are smaller and computationally more efficient.

The optimal convergence order $k+1$ is achieved for all the unknowns, as proved for the classical HDG equal-order approximation of the velocity-pressure formulation and for the more involved discretization of the Cauchy formulation based on the $\bm{M}$-decomposition.
The novelty and main advantage of the present approach relies on being able to exploit the same degree of approximation for both primal and mixed variables, in presence of the symmetric part of the gradient.
In addition, a velocity field superconverging with order $k+2$ is obtained via a local post-process procedure, exploiting the optimal convergence of the mixed and hybrid variables.
%

Numerical studies show the optimal convergence and superconvergence properties of the method in 2D and 3D using meshes of different element types and the robustness of the approach with respect to the choice of the HDG stabilization parameter.
Eventually, the drag force on a sphere is evaluated using different degrees of approximations to show the capability of the method to compute industrially relevant quantities of interest with an acceptable precision.

\section*{Acknowledgements}

This work was partially supported by the European Union's Horizon 2020 research and innovation programme under the Marie Sk\l odowska--Curie grant agreement No. 675919 and the Spanish Ministry of Economy and Competitiveness (Grant number: DPI2017-85139-C2-2-R).
The support of the Generalitat de Catalunya (Grant number: 2017SGR1278) is also gratefully acknowledged.
Finally, Alexandros Karkoulias was supported by the European Education, Audiovisual and Culture Executive Agency (EACEA) under the Erasmus Mundus Joint Doctorate Simulation in Engineering and Entrepreneurship Development (SEED), FPA 2013-0043.

\bibliographystyle{plain}
\bibliography{Ref-HDG}

\end{document}